\documentclass{amsart}

\usepackage{graphicx}
\usepackage{amscd}
\usepackage{epstopdf}
\usepackage{amssymb}
\theoremstyle{plain}

\newtheorem{thm}{Theorem}

\newtheorem{cor}[thm]{Corollary}
\newtheorem{lem}[thm]{Lemma}

\theoremstyle{definition}
\newtheorem{df}[thm]{Definition}
\newtheorem{rem}[thm]{Remark}
\newtheorem{eg}[thm]{Example}
\newtheorem{calc}[thm]{Calculation}

\author{Michael Robinson} 
\address{Center for Applied Mathematics\\657 Rhodes Hall\\Cornell
  University, Ithaca, NY 14850}
\email{robinm@cam.cornell.edu}


\title[Asymptotic-numerical approach for ODE]{An asymptotic-numerical
  approach for examining global solutions to an ordinary differential
  equation}

\begin{document}

\begin{abstract}
  Purely numerical methods do not always provide an accurate way to
  find all the global solutions to nonlinear ODE on infinite
  intervals.  For example, finite-difference methods fail to capture
  the asymptotic behavior of solutions, which might be critical for
  ensuring global existence.  We first show, by way of a detailed
  example, how asymptotic information alone provides significant
  insight into the structure of global solutions to a nonlinear ODE.
  Then we propose a method for providing this missing asymptotic data
  to a numerical solver, and show how the combined approach provides
  more detailed results than either method alone.
\end{abstract}

\maketitle

\section{Introduction}

Finding global solutions to nonlinear ordinary differential equations
on an infinite interval can be rather difficult.  Numerical
approximations can be particularly misleading, especially because they
examine only a finite-dimensional portion of the infinite-dimensional
space in which solutions lie.  Additionally, the conditions for global
existence can be rather delicate, which a numerical solver may have
difficulty rigorously checking.  In situations where there is
well-defined asymptotic behavior for global solutions, it is possible
to exploit the asymptotic information to answer questions about global
existence and uniqueness of solutions directly.  Additionally, more
detailed information may be provided by using the asymptotic behavior
to install artificial boundary conditions for use in a numerical
solver.  The numerical solver can then run on a bounded interval with
boundary conditions that match the numerical approximations to an
asymptotic expansion valid on the rest of the solution interval.

For concreteness, we consider the behavior of solutions satisfying the
differential equation
\begin{equation}
\label{long_ode1}
0=f''(x)-f^2(x)+\phi(x),\text{ for all }x\in\mathbb{R}.
\end{equation}
In particular, we wish to know how many solutions there are for a
given $\phi$.  (There may be uncountably many solutions, as in the
case where $\phi \equiv \text{const}>0$.)  This problem depends rather
strongly on the asymptotic behavior of solutions to \eqref{long_ode1}
as $|x| \to \infty$, so it is useful to study instead the pair of
initial value problems
\begin{equation}
\label{ode1}
\begin{cases}
0=f''(x) - f^2(x) + \phi(x) \text{ for }x>0\\
(f(0),f'(0))\in Z,
\end{cases}
\end{equation}
and
\begin{equation}
\begin{cases}
\label{ode1_backwards}
0=f''(x) - f^2(x) + \phi(x) \text{ for }x<0\\
(f(0),f'(0))\in Z',
\end{cases}
\end{equation}
where $\phi \in C^\infty(\mathbb{R}).$ The sets $Z$,$Z'$ supply the
  initial conditions for which solutions exist to \eqref{ode1} for all
  $x>0$ and to \eqref{ode1_backwards} for all $x<0$, respectively.
  Solutions to \eqref{long_ode1} will occur exactly when $Z \cap Z'$
  is nonempty.  Indeed, the theorem on existence and uniqueness for
  ODE gives a bijection between points in $Z\cap Z'$ and solutions to
  \eqref{long_ode1}. \cite{LeeSmooth} Since \eqref{ode1} and
  \eqref{ode1_backwards} are related by reflection across $x=0$, it is
  sufficient to study \eqref{ode1} only.

Due to the asympotic behavior of solutions to \eqref{ode1}, the
methods we employ here will be most effective in the specific cases
where $\phi$ is nonnegative and montonically decreasing to zero.  (We
denote the space of smooth functions that decay to zero as
$C^\infty_0(\mathbb{R})$.)  The decay condition on $\phi$ allows the
differential operators in \eqref{long_ode1} through
\eqref{ode1_backwards} to be examined with a perturbative approach as
$x$ becomes large, and makes sense if one is looking for smooth
solutions in $L^p(\mathbb{R})$ with bounded derivatives.  

When $\phi$ is strictly negative, it happens that no solutions exist
to \eqref{ode1} for all $x>0$.  The monotonicity restriction on $\phi$
provides some technical simplifications and sharpens the results that
we obtain.  This leads us to restrict $\phi$ to a class of functions
that captures this monotonicity restriction but allows some
flexibility, which we shall call the {\it M-shaped functions}.

It is unlikely that we will be able to solve \eqref{ode1} explicitly
for arbitrary $\phi$, so one might think that numerical approximations
might be helpful.  However, most numerical approximations will not be
able to count the number of global solutions accurately.  For
instance, finite-difference methods are typically only useful for
finding solutions valid on finite intervals of $\mathbb{R}$.  This is
unfortunately not sufficient, since the behavior of solutions to
\eqref{ode1} will be shown in Theorem \ref{limzero_iff_bounded_lem} to
either tend to zero or fail to exist.  A typical finite-difference
solution that {\it appears} to tend to zero may in fact not, and as a
result fail to be a solution over all $x>0$.  

Because of this failure, we need to understand the asymptotic behavior
of solutions to \eqref{ode1} as we take $x \to \infty$.  Equivalently,
since $\phi \in C_0^\infty(\mathbb{R})$, this means that we should
examine solutions with $\phi$ small.  The driving motivation for this
discussion is that solutions to $0=f''(x)-f^2(x)+\phi(x)$ for $\phi$
small behave much like solutions to $0=f''(x)-f^2(x)$.  In the latter
case, we can completely characterize the solutions which exist on
intervals like $[x_0,\infty)$.

In Section \ref{review_sec} we review what is known about the much
simpler case where $\phi$ is a constant.  Of course, then \eqref{ode1}
is autonomous, and the results are standard.  In Section
\ref{asymp_exist_sec}, we establish the existence of solutions which
are asymptotic to zero.  Some of these solutions are computed
explicitly using perturbation methods in Section \ref{series_sec},
where low order approximations are used to gather qualitative
information about the initial condition sets $Z$ and $Z'$.  In
Sections \ref{extension_sec} and \ref{geom_props_Z}, these qualitative
observations are made precise.  Section \ref{long_sec} applies these
observations about $Z$ and $Z'$ to give existence and uniqueness
results for \eqref{long_ode1}.  Finally, in Section \ref{numer_sec},
we use the information gathered about $Z$ and $Z'$ to provide
artificial boundary conditions to a numerical solver on a bounded
interval, which sharpens the results from Section \ref{long_sec}.  We
exhibit the numerical results for a typical family of $\phi$, showing
bifurcations in the global solutions to \eqref{long_ode1}.

\section{Review of behavior of solutions to $0=f''(x)-f^2(x)+P$}
\label{review_sec}

It will be helpful to review the behavior of 
\begin{equation}
\label{IVP2}
\begin{cases}
0=f''(x) - f^2(x) + P\\
f(0),f'(0)\text{ given},
\end{cases}
\end{equation}
where $P$ is a constant, since varying $\phi$ can be viewed as a
perturbation on the case $\phi(x)=P$.  In particular, we need to
compute some estimates for later use.  We shall typically take $P>0$,
as there do not exist solutions for all $x$ if $P<0$.  




\begin{figure}
\includegraphics[height=3in]{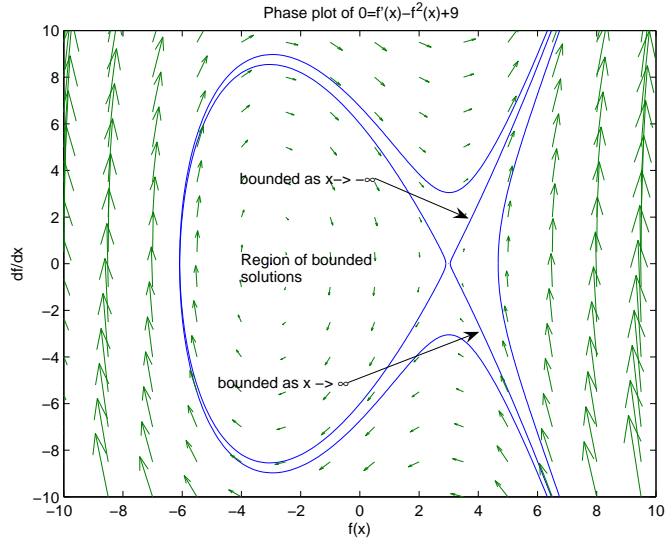}
\caption{The phase plot of $f''-f^2+9=0$.  Bounded solutions live in a
small region, the rest are unbounded.}
\label{P_eq_p9}
\end{figure}

\begin{lem}
\label{no_upper_f_bnd_lem}
Suppose $f$ is a solution to the initial value problem
\eqref{IVP2} with $f(0)>\sqrt{P}$ and $f'(0)>0$.  Then there does
not exist an upper bound on $f(x)$, when $x>0$.  Additionally, if
$P<0$, there does not exist an upper bound on $f(x)$.
\end{lem}

\begin{df}
The differential equation \eqref{IVP2} comes from a Hamiltonian, namely
\begin{equation*}
H(f,f')=\frac{1}{3} f^3 - \frac{1}{2}f'^2-fP+\frac{2}{3}P^{3/2}.
\end{equation*}
\end{df}

\begin{lem}
\label{bounded_in_funnel_lem}
Suppose $f$ is a solution to the equation \eqref{IVP2} on
$\mathbb{R}$.  All bounded solutions lie in the funnel
\begin{equation}
\label{M}
M=\{(f,f')|H(f,f')\ge 0 \text{ and } f \le \sqrt{P}\}.
\end{equation}
Any solution which includes a point outside the closure of $M$ is
unbounded, either for $x>0$ or $x<0$.  (Note that $M$ is the
teardrop-shaped region in Figure \ref{P_eq_p9}.)
\begin{proof}
\begin{itemize}
\item $M$ is a bounded set.  Notice that $H(f,0) \ge H(f,f')$, or in other
  words within $M$,
\begin{equation*}
0<\frac{1}{3}f^3 - \frac{1}{2}f'^2-f P + \frac{2}{3}P^{3/2}\le
 \frac{1}{3} f^3 - fP + \frac{2}{3}P^{3/2}.
\end{equation*}
Elementary calculus reveals that this inequality establishes a
lower bound on $f$, namely that 
\begin{equation}
\label{f_bound}
-\sqrt{3P} \le f \le \sqrt{P}
\end{equation} 
On the other hand,
\begin{equation}
\label{fp_bound}
|f'|<\sqrt{\frac{4}{3}P^{3/2}+\frac{2}{3}f^3-2fP} \le \sqrt{\frac{8}{3}}P^{3/4}
\end{equation}
immediately establishes a bound on $f'$.

\item $M$ is a funnel, from which solutions neither enter nor leave.
This is immediate from the fact that $H$ is the Hamiltonian, and the
definition of $M$ simply says that $H(f,f')>0$.
This suffices since solutions to \eqref{IVP2} are tangent to level
curves of $H$.

\item If $(f(0),f'(0)) \notin M$ then $f$ is unbounded.  Evidently if
  $f(0)> \sqrt{P}$ and $f'(0)>0$, then Lemma \ref{no_upper_f_bnd_lem}
  applies to give that $f$ is unbounded.  For the remainder, discern
  two cases.  First, suppose $f(0)>\sqrt{P}$ and $f'(0)<0$.
  Evidently, $H(f(0),f'(0))=H(f(0),-f'(0))$, so it's just a matter of
  verifying that a solution curve transports our solution to the first
  quadrant.  But this is immediately clear from the formula for
\begin{equation*}
f'=\pm \sqrt{\frac{2}{3}f^3-2fP-2H(f(0),f'(0))},
\end{equation*}
which gives $f'=\pm f'(0)$ when $f=f(0)$.  The other case is when
   $H(f(0),f'(0)) \le 0$.  Then we show that there
   is a point $(\sqrt{P},g)$ on the same solution curve, and then
   Lemma \ref{no_upper_f_bnd_lem} applies.  So we try to satisfy
\begin{eqnarray*}
\frac{1}{3}P^{3/2} - \frac{1}{2}g^2-P^{3/2} + \frac{2}{3} P^{3/2} &=&
H(f(0),f'(0)) \le 0 \\
g^2 &=& -2 H(f(0),f'(0)) \ge 0,
\end{eqnarray*}
which clearly has a solution in $g$.  Finally, if $g=0$, then
$f(0)>\sqrt{P}$, so it has already been covered above.
\end{itemize}
\end{proof}
\end{lem}

\begin{lem}
\label{asymptote_lem}
If $f$ is a solution to \eqref{IVP2} with $f(0)>\sqrt{P}$, and
$f'(0)>0$ then there exists a $C$ such that $\lim_{x \to C} f(x)
= \infty$. 
\end{lem}

\section{Existence of asymptotic solutions for $\phi \in
  C_0^\infty(\mathbb{R})$} 
\label{asymp_exist_sec}

The first collection of results we obtain will make the assumption
that $\phi$ tends to zero.  From this, a number of useful asymptotic
results follow.  Working in the phase plane will be useful for
understanding \eqref{ode1}.  Of course \eqref{ode1} is not autonomous,
but by adding an additional variable, it becomes so.

\begin{df}
We think of \eqref{ode1} as a vector field $V$ on $\mathbb{R}^3$, defined
by the formula
\begin{equation}
\label{ode1_system}
V(f,f',x)=\begin{pmatrix}
f' \\ f^2 - \phi(x) \\ 1
\end{pmatrix}.
\end{equation}
Notice that the first coordinate of an integral curve for this vector
field solves \eqref{ode1}.
\end{df}

\begin{df}
Define $H(f,f',x)=\frac{1}{3}f^3-\frac{1}{2}f'^2-f \phi(x) +
\frac{2}{3} \phi^{3/2}(x).$  Notice that for constant $\phi = P$, this
reduces to a Hamiltonian for \eqref{IVP2}.  
\end{df}

\begin{thm}
\label{limzero_iff_bounded_lem}
Suppose $f$ is a solution to the problem \eqref{ode1} where $\phi \in
C_0^\infty(\mathbb{R})$.  If $f$ does not tend to zero as
$x\to \infty$, then there exists a $z$ such that $\lim_{x\to z} f(x) =
  \infty$.  Stated another way, if $f$ solves \eqref{ode1} for all
  $x>0$, then $\lim_{x \to \infty} f(x) = 0$.
\begin{proof}
If $f$ does not tend to zero, this means that there is an $R>0$ such
that for each $x_0>0$, there is an $x>x_0$ so that $|f(x)|>R$.  But
since $\phi$ tends to zero as $x \to \infty$, for any $P>0$ we can
find an $x_1>0$ such that for all $x>x_1$, $|\phi(x)|<P$.  Choose such
a $P$ so that the set $M$ in Lemma \ref{bounded_in_funnel_lem}
associated to \eqref{IVP2} is contained entirely within the strip
$-R<f<R$.  We can do this since the set $M$ is bounded, and its radius
decreases with decreasing $P$, as shown in \eqref{f_bound} and
\eqref{fp_bound}.  But this means that there is an $x_2>x_1$ such that
$|f(x_2)|>R$.

Construct the following regions (See Figure \ref{limzero_fig}):
\begin{equation*}
I=\{(f,f',x)| f \ge R \text{ and } f' \le 0\},
\end{equation*}
\begin{equation*}
II=\{(f,f',x)| f \ge R \text{ and } f' \ge 0\},
\end{equation*}
\begin{equation*}
III=\{(f,f',x)| f \le R \},
\end{equation*}
and
\begin{equation*}
IV=\left\{(f,f',x)| f' \ge 0 \text{ and } f'\ge -R \text{ and } \left(
\frac{1}{3}f^3 - \frac{1}{2} f'^2 - f P + \frac{2}{3} P^{3/2} \le 0
\text{ if } f \le \sqrt{P}\right) \right\}.
\end{equation*}

\begin{figure}
\includegraphics[height=3in]{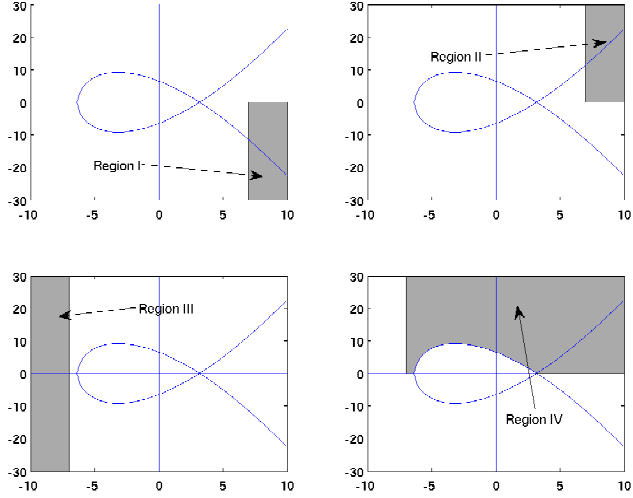}
\caption{The Regions $I$, $II$, $III$, and $IV$ of Theorem \ref{limzero_iff_bounded_lem}}
\label{limzero_fig}
\end{figure}

The following statements hold:
\begin{itemize}
\item Region $I$ is an antifunnel.  Along $f=R$ and $f'=0$, solutions
  must exit.  Once a solution exits Region $I$, it cannot reenter.
  Also, because $f > \sqrt{P}$, $f''=f^2-\phi> f^2 - P > 0$, solutions
  must exit Region $I$ in finite $x$.
\item Region $II$ is a funnel.  Along $f=R$ and $f'=0$, solutions
  enter.  Now $f''=f^2 - \phi > f^2 - P \ge 0$ and $f' \ge 0$, so
  solutions will increase at an increasing rate and so, they
  are unbounded. 
\item Solutions remain in Region $III$ for only finite $x$, after
  which they must enter Region $IV$.  This occurs since $f \le
  -\sqrt{P} < 0$, and so $f'$ always increases.  Note that for $f'<0$,
  solutions will enter Region $III$ along $f=-R$, and for $f'>0$,
  solutions exit along $f=-R$.
\item Region $IV$ is a funnel.  Solutions enter along $f=-R$ and along
  $f'=0$ (note that $|f| \ge \sqrt{P}$ in both cases).  Along the
  curve boundary of Region $IV$, we have that
\begin{eqnarray*}
\nabla \left(\frac{1}{3} f^3 - \frac{1}{2} f'^2 - fP +
\frac{2}{3}P^{3/2} \right) \cdot  V(f,f',x)&=&
\begin{pmatrix} f^2 - P \\ -f' \\ 0 \end{pmatrix}^T \begin{pmatrix} f' \\ f^2 - \phi \\ 1
\end{pmatrix}\\
&=&f'(P-\phi)< 0,
\end{eqnarray*}
so that solutions enter.
\end{itemize}

Now suppose $(f(x_2),f'(x_2),x_2)\in I$.  After finite $x$, say at
$x=x_2'$, the solution through that point must exit Region $I$, never
to return.  Then, there is an $x_3 > x_2'$ such that $|f(x_3)|>R$.  So
this solution has either $(f(x_3),f'(x_3),x_3) \in II$ or $\in III$.
The former gives the conclusion we want, so consider the latter case.
The solution will only remain in Region $III$ for finite $x$, after
which it enters Region $IV$, say at $x=x_3'$.  Then there is an
$x_4>x_3'$ such that $|f(x_4)|>R$.  Now the only possible location for
$(f(x_4),f'(x_4),x_4)$ to be is within Region $II$, since it must also
remain in Region $IV$.  As a result, the solution is unbounded by an
easy extension of Lemma \ref{no_upper_f_bnd_lem}.  As $x$ becomes
large, $\phi$ tends to zero, so the solution will be asymptotic to an
unbounded solution of $0=f''-f^2$.  But Lemma \ref{no_upper_f_bnd_lem}
above assures us that such a solution is unbounded from above, and
Lemma \ref{asymptote_lem} gives that it has an asymptote.  Hence, our
solution must blow up at a finite $x$.
\end{proof}
\end{thm}

This result indicates that solutions to \eqref{ode1} which exist for
all $x>0$ are rather rare.  Those which exist for all $x>0$ must tend
to zero, and it seems difficult to ``pin them down.''  We now apply
topological methods, similar to those employed in \cite{HubbardWest},
to ``capture'' the solutions we seek.

We begin by extending the usual definition of a flow slightly to the
case of a manifold with boundary.

\begin{df}
Suppose $M$ is a manifold with boundary.  A {\it flow domain} $J$ is a
  subset of $\mathbb{R} \times M$ such that if $x \in J$ then
  $J_x=\text{pr}_1 (J \cap \mathbb{R} \times \{x\})$ is an interval
  containing 0, and if $x$ is in the interior of $M$ then 0 is
  in the interior of $J_x$.  ($\text{pr}_1:\mathbb{R}\times M \to
  \mathbb{R}$ is projection onto the first factor)
\end{df}

\begin{df}
A {\it (smooth) flow} is a smooth map $\Phi$ from a flow domain $J$
to a manifold with boundary $M$, satisfying
\begin{itemize}
\item $\Phi(0,x)=x$ for all $x \in M$ and
\item $\Phi(t_1+t_2,x)=\Phi(t_1,\Phi(t_2,x))$ whenever both
  sides are well-defined.
\end{itemize}
Additionally, we assume that flows are {\it maximal} in the sense that
they cannot be written as a restriction of a map from a larger flow
domain which satisfies the above axioms.  We call the curve $\Phi_x:J_x
\to M$ defined by $\Phi_x(t)=\Phi(t,x)$ the {\it integral curve
through $x$} for $\Phi$.
\end{df}

\begin{df}
Suppose $\Phi:J \to M$ is a flow on $M$ and $x \in \partial M$.  Then
the flow at $x$ is said to be {\it inward-going} (or simply {\it
inward}) if $J_x$ is an interval of the form $[0,a)$ or $[0,a]$ for
some $0<a\le \infty$.  Likewise, the flow at $x$ is {\it outward-going} if
$J_x$ is of the form $(a,0]$ or $[a,0]$ for $-\infty \le a < 0$.
\end{df}

\begin{thm} (Antifunnel theorem)
Suppose $\Phi:J \to M$ is a flow on $M$ and that $\{A,B\}$ forms a
partition of the boundary of $M$ such that the flow of $\Phi$ is
inward along A and outward along B.  If every integral curve of $\Phi$
intersects $B$ in finite time (ie. $J_x$ is bounded for each $x$),
then $A$ is diffeomorphic to $B$.
\begin{proof}
For each $x \in A$, $J_x=[0,t_x]$, where $t_x$ is the time
which the integral curve through $x$ intersects $B$.  (We have that
$\Phi(t_x,x)$ is outward-going, since $J_x$ is closed, so it is in
$B$.)

Using this, we can define a map $F:A \to B$ by $F(x)=\Phi(t_x,x)$.
Claim that $F$ takes $A$ smoothly and injectively into $B$.  The
smoothness follows from the smoothness of $\Phi$ and that $\partial M$
is a smooth submanifold.  To see the injectivity, suppose $F(x)=F(y)$ for
some $x,y \in A$, so $\Phi(t_x,x)=\Phi(t_y,y)$.  Without loss of
generality, suppose $0 < t_x \le t_y$.  Then we have that
\begin{eqnarray*}
F(x)&=&F(y)\\
\Phi(-t_x,F(x))&=&\Phi(-t_x,F(y))\\
\Phi(-t_x,\Phi(t_x,x))&=&\Phi(-t_x,\Phi(t_y,y))\\
\Phi(t_x-t_x,x)&=&\Phi(t_y-t_x,y)\\
x&=&\Phi(t_y-t_x,y).
\end{eqnarray*}
But the flow is inward at $x$, so it is also inward at
$\Phi(t_y-t_x,y)$.  This means that $(t_y - t_x - \epsilon,y) \notin J$
for every $\epsilon > 0$.  But this contradicts the fact that $(t_y,y)\in J$
unless we have $t_y \le t_x$.  As a result, $t_y = t_x$, so $x=y$.

In just the same way as for $F$, we construct a map $G:B \to A$ so
that $G$ takes $B$ smoothly and injectively into $A$.  Namely, we
suppose $J_y=[s_y,0]$ for some $s_y$, and put $G(y)=\Phi(s_y, y)$.
Notice that by maximality, if there were to be an $x \in A$ such that
$F(x)=y$, $s_y = - t_x$.

Now we claim that $G$ is the inverse of $F$. We have that 
\begin{eqnarray*}
(G\circ F)(x)&=&\Phi(s_{F(x)},F(x)) \\
&=&\Phi(s_{F(x)},\Phi(t_x,x))\\
&=&\Phi(s_{F(x)}+t_x,x)\\
&=&\Phi(-t_x+t_x,x)=x,\\
\end{eqnarray*}
where we employ the remark about $s_y$ above.
\end{proof}
\end{thm}

\begin{rem}
We can extend the Antifunnel theorem to a topological space $X$ on
which a flow $\Phi:J \to X$ acts in the obvious way.  In that case,
there is no reasonable definition of the boundary of $X$.  However,
the notion of inward- and outward-going points still makes sense.  If
we let $A$ be the set of inward-going points and $B$ be the set of
outward-going points in $X$, then the conclusion is that $A$ is
homeomorphic to $B$.
\end{rem}

Now we employ the Antifunnel theorem to deduce the existence of
a bounded solution to $0=f''-f^2+\phi$ for $x>x_0$ for some $x_0\ge 0$.

\begin{thm}
\label{region_r1_thm}
Suppose $0 \le \phi(x) \le K$ for all $x \ge x_0$ for some $x_0>0$ and
$0 < K < \infty$.  Then the region $R_1$ given by
$R_1=\{(f,f',x)|H(f,f',x) \ge 0, x \ge 0, f\le \sqrt{\phi(x)}\}$
contains a bounded solution to $0=f''(x)-f^2(x)+\phi(x)$, which exists
for all $x$ greater than some nonnegative $x_1$.
\begin{proof}
Without loss of generality, we may take $x_0=0$, because otherwise
solutions must exit the portions of $R_1$ in $\{(f,f',x)|x<x_0\}$
since the $x$-component of $V(f,f',x)$ is equal to 1.

If $\phi(0)>0$, partition the boundary of $R_1$ into two pieces:
$A=\{(f,f',x)|x = 0\}$ and $B=\{(f,f',x)|H(f,f',x)=0\}$.  The flow of
$V$ is evidently inward along $A$.  As for $B$, notice that $\nabla H$
is an inward-pointing vector field normal to $B$.  We compute
\begin{eqnarray*}
\nabla H \cdot V &=& 
\begin{pmatrix} f^2 - \phi(x) \\ -f' \\ (-f + \sqrt{\phi(x)})\phi'(x) 
  \end{pmatrix}^T
\begin{pmatrix} f' \\ f^2 - \phi(x) \\ 1 \end{pmatrix}\\
&=& ( - f + \sqrt{\phi(x)} ) \phi'(x),
\end{eqnarray*}
which has the same sign as $\phi'(x)$ when $f < \sqrt{\phi(x)}$ in
$R_1$.  Finally, we must deal with the case where $f=\sqrt{\phi(x)}
\in B$.  But in this case, $f'=0$ from the equation for $H$, so we see
that $V(\sqrt{\phi(x)},0,x)=(0,0,1)^T$, so the flow is inward when
$\phi'(x)<0$ and outward when $\phi'(x)>0$.  This means that
the portion of the boundary of $R_1$ on which the flow is outward is a
disjoint union of annuli.  On the other hand, the portion of the
boundary of $R_1$ on which the flow is inward is the disjoint union of
a disk (namely $R_1 \cap \{x=0\}$) and some annuli.  Now if
$\phi(0)\le 0$, then the set $A$ above is empty.  However, it is
obvious that set the inflow portion of the boundary is still
homeomorphic to the disjoint union of a disk and some annuli.

We can apply the Antifunnel theorem to conclude that there is a
solution which does not intersect either the inflow or outflow
portions of the boundary.  There is a lower bound on the
$x$-coordinate of such a solution, since the $x$-component of
$V(f,f',x)$ is equal to 1, and the Region $R_1$ lies within the
half-space $x>0$.  Therefore, there must exist a solution which
enters $R_1$, and remains inside the interior of $R_1$ for all larger
$x$.  That such a solution is bounded follows from the fact that each
constant $x$ cross section of $R_1$ has a radius bounded by the
inequalities \eqref{f_bound} and \eqref{fp_bound}, and the fact that
$\phi(x) \le K < \infty$.  
\end{proof}
\end{thm}

\section{Asymptotic series solution}
\label{series_sec}

Theorem \ref{region_r1_thm} ensures the existence of solutions to
$0=f''-f^2+\phi$ for $x$ sufficiently large.  However, it does not
give any description of the initial condition set $Z$ which leads to
such solutions, nor does it give a description of the maximal
intervals of existence.  Fortunately, it is relatively easy to
construct an asymptotic series for solutions to \eqref{ode1}, which
will provide a partial answer to this concern.  In doing so, we
essentially follow standard procedure, as outlined in \cite{Holmes},
for example.  However, our case is better than the standard situation,
because under relatively mild restrictions this series {\it converges}
to a true solution.

We begin by supposing that our solution has the form
\begin{equation}
\label{f_format}
f=\sum_{k=0}^\infty f_k,
\end{equation}
where we temporarily assume $f_{k+1} \ll f_k$ and $f_0 \gg \phi$, as
$x \to +\infty$.  (This assumption will be verified in Lemma
\ref{series_bound_lem}.)  Substituting \eqref{f_format} into
\eqref{ode1}, we get

\begin{eqnarray*}
0&=&\sum_{k=0}^\infty \left[ f_{k}'' - \sum_{m=0}^{k}f_m
  f_{k-1}\right]+ \phi \\
0&=&f_0''-f_0^2+(f_1''-2f_0f_1+\phi)+\sum_{k=2}^\infty
  \left[f_{k}''-2f_0f_{k}-\sum_{m=1}^{k-1}f_mf_{k-m}\right].
\end{eqnarray*}
We solve this equation by setting different orders to zero.
Namely,
\begin{eqnarray*}
0&=&f_0''-f_0^2\\
0&=&f_1''-2f_0f_1+\phi\\
0&=&f_k''-2f_0f_k-\sum_{m=1}^{k-1}f_mf_{k-m}.
\end{eqnarray*}
The equation for $f_0$ is integrable, and therefore easy to solve.
(There are two families of solutions for $f_0$.  We select the
nontrivial one, because the other one simply results in $f(x) \sim
-\int_x^\infty \int_t^\infty \phi(s) ds\, dt$.)  The equations for
$f_k$ are linear and can be solved by a reduction of order.  Thus
formally, the solutions are
\begin{equation}
\label{series_soln_coef}
\begin{cases}
f_0=\frac{6}{(x-d)^2}\\
f_1=\frac{1}{(x-d)^3}\left[ K + \int^x(t-d)^6\int_t^\infty
  \frac{\phi(s)}{(s-d)^3}ds\,dt\right]\\
f_k=-\frac{1}{(x-d)^3}\int^x (t-d)^6 \int_t^\infty \frac{ \sum_{m=1}^{k-1}
  f_m(s) f_{k-m}(s)}{(s-d)^3}ds\,dt,
\end{cases}
\end{equation}
for $d,K$ constants.  Notice that these constants parametrize the set of
initial conditions $Z$.  

\begin{lem}
\label{series_bound_lem}
Suppose $f(x)=\sum_{k=0}^\infty f_k(x)$ where the $f_k$ are given by
\eqref{series_soln_coef}.  If there
exists an $M>0$, an $R>0$, and an $\alpha>5$ such that 
\begin{equation}
\label{phi_cond}
|\phi(x)|<\frac{M}{(x-d)^\alpha} \text{ for all } |x-d|>R>0,
\end{equation}
then $f(x)$ is bounded above by the power series
\begin{equation}
\label{f_power_series_bound}
|f(x)| \le \frac{1}{(x-d)^2} \sum_{k=0}^\infty \left|\frac{A_k}{x-d}\right|^k.
\end{equation}
\begin{proof}
We proceed by induction, and begin by showing that the $f_1$ term is
appropriately bounded:

\begin{eqnarray*}
|f_1(x)|&=&\left|\frac{1}{(x-d)^3}\left[K+\int^x(t-d)^6\int_t^\infty
  \frac{\phi(s)}{(s-d)^3}ds\,dt\right]\right|\\
&\le&\left|\frac{1}{(x-d)^3}\left[K+\int^x(t-d)^6\int_t^\infty
  \frac{M}{(s-d)^{3+\alpha}}ds\,dt\right]\right|\\
&\le&\left|\frac{1}{(x-d)^3}\left[K+
  \frac{M}{(2+\alpha)(5-\alpha)(x-d)^{\alpha-5}} \right]\right|\\
\end{eqnarray*}
Now since $|x-d|>R$ and $\alpha > 5$, we have that
\begin{eqnarray*}
|f_1(x)|&\le& \frac{1}{|x-d|^3} \left[|K|+
  \frac{M}{(2+\alpha)|5-\alpha|R^{\alpha-5}} \right].\\
&\le&\frac{A_1}{|x-d|^3}.
\end{eqnarray*}
with
\begin{equation}
\label{A_1_def}
A_1 = |K|+\frac{M}{(2+\alpha)(\alpha-5) R^{\alpha-5}}.
\end{equation}

For the induction hypothesis, we assume that
$|f_i|\le\frac{A_i}{|x-d|^{2+i}}$ with $A_i \ge 0$ and for all $i \le
k-1$.  We have that

\begin{eqnarray*}
\sum_{m=1}^{k-1} f_m f_{k-m} &\le& \sum_{m=1}^{k-1}
\frac{A_m}{|x-d|^{2+m}} \frac{A_{k-m}}{|x-d|^{2+k-m}} \\
&\le& \frac{1}{|x-d|^{k+4}} \sum_{m=1}^{k-1} A_m A_{k-m}, \\
\end{eqnarray*}
so by the same calculation as for $f_1$, we obtain
\begin{equation*}
f_k \le \frac{\sum_{m=1}^{k-1} A_m A_{k-m}}{(k+6)(k-1)}\frac{1}{|x-d|^{k+2}}.
\end{equation*}

Hence we should take
\begin{equation}
\label{A_k_def}
A_k = \frac{\sum_{m=1}^{k-1} A_m A_{k-m}}{(k+6)(k-1)}.
\end{equation}

Hence we have that 

\begin{equation*}
|f(x)| \le \sum_{k=0}^\infty |f_k(x)|\le
 \frac{1}{|x-d|^2}\sum_{k=0}^\infty A_k \left|\frac{1}{x-d}\right|^k.
\end{equation*}
\end{proof}
\end{lem}

\begin{lem}
\label{series_conv_lem}
The power series given by 
\begin{equation*}
\sum_{k=0}^\infty \frac{A_k}{|x-d|^k}, 
\end{equation*}
with $A_0, A_1 \ge 0$ given, and 
\begin{equation*}
A_{k}=\frac{\sum_{m=1}^{k-1} A_m A_{k-m}}{(k+6)(k-1)} =
\frac{\sum_{m=1}^{k-1} A_m A_{k-m}}{k^2+5k-6}
\end{equation*}
converges for $|x-d|>R$ if $A_1 \le 8R$.
\begin{proof}
We show that under the conditions given, the series passes the usual
ratio test.  That is, we wish to show
\begin{equation*}
\lim_{k \rightarrow \infty} \left|\frac{A_{k+1}}{A_k}\right| \le R.
\end{equation*}
Proceed by induction.  Take as the base case, $k=1$: by the formula for
$A_{k}$,
\begin{equation*}
A_2=\frac{A_1^2}{8}, \text{ so }
\frac{A_2}{A_1}=\frac{A_1^2}{8A_1}=\frac{A_1}{8} \le R.
\end{equation*}

Then for the induction step, 
\begin{eqnarray*}
\frac{A_{k+1}}{A_k}&=&\frac{\sum_{m=1}^k A_mA_{k-m+1}}{A_k (k^2+7k)}\\
&=&\frac{\sum_{m=2}^k A_mA_{k-m+1}+A_1 A_k}{A_k (k^2+7k)}\\
&=&\frac{\sum_{m=1}^{k-1} A_{m+1}A_{k-m}+A_1 A_k}{A_k (k^2+7k)}\\
&\le&\frac{R\sum_{m=1}^{k-1} A_mA_{k-m}+A_1 A_k}{A_k (k^2+7k)}\\
&\le&\frac{RA_k(k^2+5k-6)+A_1 A_k}{A_k (k^2+7k)}\\
&\le&\frac{R(k^2+5k-6)+A_1}{(k^2+7k)}\\
&\le&\frac{R(k^2+5k-2)}{(k^2+7k)}\\
&\le&R,\\
\end{eqnarray*}
since $A_1 \le 8 R$.  Thus $\left| \frac{A_{k+1}}{A_k} \right| \le R$
for all $k$, so the power series converges.
\end{proof}
\end{lem}

Lemma \ref{series_conv_lem} provides conditions for the convergence of the
bounding series found in Lemma \ref{series_bound_lem}.  Hence we have actually proven the
following:

\begin{thm}
\label{series_conv_thm}
Suppose $f(x)=\sum_{k=0}^\infty f_k(x)$ where the $f_k$ are given by
\eqref{series_soln_coef}.  If there exists an $M>0$, an $R>0$, an
$\alpha>5$ such that \eqref{phi_cond} holds, and furthermore
\begin{equation}
\label{series_conv_cond}
M < 8 (\alpha+2)(\alpha-5)R^{\alpha-4},
\end{equation}
then the series for $f(x)$ converges for all $x$ such that $|x-d|>R$.
\begin{proof}
Combining Lemmas \ref{series_bound_lem} and \ref{series_conv_lem}, we
find that the key condition is that $A_1 \le 8R$, which by
substitution into \eqref{A_1_def} yields
\begin{equation*}
0<|K|+\frac{M}{(\alpha+2)(\alpha-5)R^{\alpha-5}}\le 8R.
\end{equation*}
But in order to have $|K| \ge 0$, this gives
\begin{equation*}
0<8R-\frac{M}{(\alpha+2)(\alpha-5)R^{\alpha-5}},
\end{equation*}
which leads immediately to the condition stated.
\end{proof}
\end{thm}

\begin{figure}
\includegraphics[height=3in]{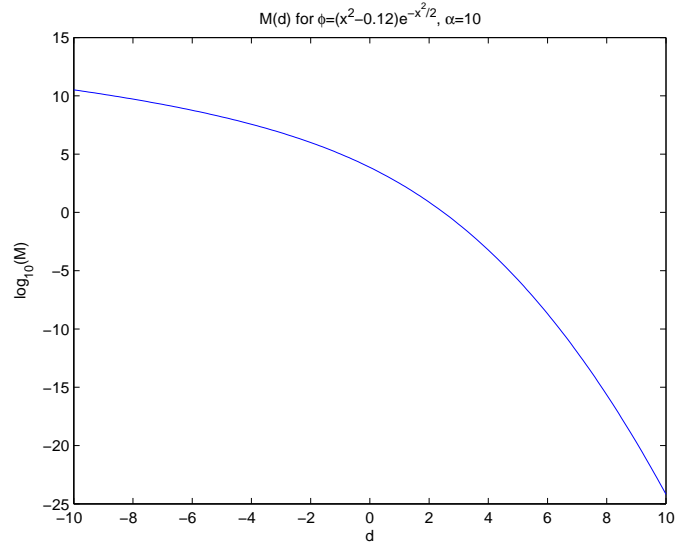}
\caption{A typical $M(d)$ function}
\label{m_of_d}
\end{figure}

\begin{figure}
\includegraphics[height=3in]{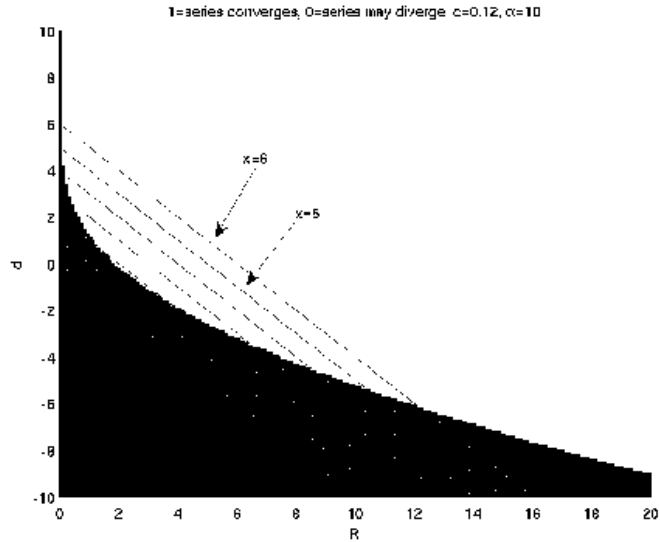}
\caption{Series convergence test, for $\phi(x)=(x^2-0.12)e^{-x^2/2}$:
  white = series converges, black = series may diverge }
\label{series_conv}
\end{figure}

\begin{eg}
It is important to notice that the $M$ defined above in Lemma
\ref{series_bound_lem} can depend crucially upon the value of $d$ and
the shape of the curve $\phi(x)$.  For the case of
$\phi(x)=(x^2-c)e^{-x^2/2}$, a typical plot of $M(d)$ is shown in Figure
\ref{m_of_d}.  It should be noted that for various values of $c$, the
$M(d)$ function is numerically very similar.

This also means that the condition \eqref{series_conv_cond} defines a
somewhat complicated region over which parameters $d, K$ and $R$ yield
convergent series solutions.  An example with our given $\phi(x)$
function is shown in Figure \ref{series_conv}.  Thus it appears that
our series solution converges if one goes out far enough, and
specifies small enough initial conditions.
\end{eg}

\begin{rem}
The convergence of the series solution is controlled by the
convergence of a well-behaved power series.  It follows that as
the $\phi$ function becomes smaller, fewer terms in the series are
needed to accurately approximate the solution.  Indeed, each term in
the series solution is asymptotically smaller than the one previous.
Thus, we can gain some qualitative information from the leading two
terms of the series, which are
\begin{equation*}
f(x) \sim \frac{6}{(x-d)^2}+\frac{1}{(x-d)^3}\left[ K + \int^x(t-d)^6\int_t^\infty
  \frac{\phi(s)}{(s-d)^3}ds\,dt\right].
\end{equation*}
Taking a derivative by $x$ gives
\begin{equation*}
f'(x) \sim \frac{-12}{(x-d)^3}+\frac{-3}{(x-d)^4}\left[ K + \int^x(t-d)^6\int_t^\infty
  \frac{\phi(s)}{(s-d)^3}ds\,dt\right]+(x-d)^3\int_x^\infty
  \frac{\phi(s)}{(s-d)^3}ds.
\end{equation*}
On the other hand, using the standard expansion for $(a+b)^{3/2}$, one
obtains
\begin{eqnarray*}
f^{3/2}(x) &\sim& \left( \frac{6^3}{(x-d)^6}+\frac{3\cdot
  6^2}{(x-d)^7}\left[ K + \int^x(t-d)^6\int_t^\infty
  \frac{\phi(s)}{(s-d)^3}ds\,dt\right] \right)^{1/2}\\
&\sim&\frac{6^{3/2}}{(x-d)^3}+\frac{(x-d)^3}{2\cdot 6^{3/2}}\frac{3\cdot
  6^2}{(x-d)^7}\left[ K + \int^x(t-d)^6\int_t^\infty
  \frac{\phi(s)}{(s-d)^3}ds\,dt\right]\\
&\sim&\frac{6^{3/2}}{(x-d)^3}+\frac{3\cdot
  18}{6^{3/2}(x-d)^4}\left[ K + \int^x(t-d)^6\int_t^\infty
  \frac{\phi(s)}{(s-d)^3}ds\,dt\right]\\
\end{eqnarray*}
which leads to 
\begin{equation}
\label{handy_asymp_exp}
f'(x) \sim -\sqrt{\frac{2}{3}} f^{3/2} + (x-d)^3\int_x^\infty
  \frac{\phi(s)}{(s-d)^3}ds.
\end{equation}
Notice that this equation depends only on $d$, not $K$.  So from this
we should expect that the initial data for solutions to be confined to
a thin region in the plane $x=0$.  This will be confirmed in Theorem
\ref{properties_of_Z}

Additionally, the relation $f'=-\sqrt{2/3} f^{3/2}$ holds exactly for
the bounded solutions of $0=f''-f^2$.  Indeed, in that case, the set
$Z$ is $\{(f,f')|3f'^2=2f^3,f'<0\}$.  So \eqref{handy_asymp_exp}
indicates that the presence of $\phi \neq 0$ will deflect the set $Z$
largely in the $f'$ direction.  This is exactly what we show in Section
\ref{geom_props_Z}.
\end{rem}

\section{Restriction to $\phi$ nonnegative and monotonically
  decreasing}
\label{extension_sec}

We now examine what stronger results can be obtained by requiring
$\phi(x)\ge 0$ and $\phi'(x)<0$ for all $x>0$.  This can be expected
to provide stronger results, in particular because the region $R_1$
employed in Theorem \ref{region_r1_thm} acquires a simpler inflow and
outflow structure on the boundary, and in particular, solutions will
exist for all $x>0$.  A collection of four results indicate that all
bounded solutions to \eqref{ode1} lie within a narrow region.

\begin{lem}
\label{region_r1_lem}
Suppose $\phi(x)\ge 0$ and $\phi'(x) < 0$ for all $x \ge 0$.  Then the
region given by $R_1=\{(f,f',x)|H(f,f',x) \ge 0, x \ge 0, f\le
\sqrt{\phi(x)}\}$ contains a bounded solution to \eqref{ode1}.
\begin{proof}
Following the proof of Theorem \ref{region_r1_thm}, we partition the
boundary of $R_1$ into two pieces: $A=\{(f,f',x)|x = 0\}$ and
$B=\{(f,f',x)|H(f,f',x)=0\}$, noting that the flow of $V$ is inward
along $A$.  Reviewing the computation in Theorem \ref{region_r1_thm}, 
the flow is outward along all of $B$.

Now we employ the Antifunnel theorem, noting that while $A$ is
simply-connected, $B$ is not.  Hence they cannot be homeomorphic, and
so there must be a solution that remains inside $R_1$ (which evidently
starts on $A$).  But the first coordinate of such an integral curve
must obviously be bounded, since the $x$ cross-sections of $R_1$ form
a decreasing sequence of sets, ordered by inclusion, and the
cross-section for $x=0$ is a bounded set.
\end{proof}
\end{lem}

\begin{lem}
\label{region_r2_lem}
Suppose $\phi(x)\ge 0$ and $\phi'(x) < 0$ for all $x \ge 0$.  Then the
region given by $R_2=\{(f,f',x)|H(f,f',x) \le 0, \frac{1}{3}f^3 -
\frac{1}{2} f'^2 \ge 0, x \ge 0,f'\le 0\}$ contains a bounded solution
to \eqref{ode1}.
\begin{proof}
Partition the boundary of $R_2$ into two pieces:
\begin{equation*}
A=\{(f,f',0)|f' \le 0\}\cup\{(f,f',x)|H(f,f',x)=0,f\le
\sqrt{\phi(x)},f'\le 0\},
\end{equation*}
and
\begin{equation*}
B=\{(f,f',x)|H(f,f',x)=0,f\ge \sqrt{\phi(x)},f'\le 0 \}\cup\{(f,f',x)|
\frac{1}{3}f^3 - \frac{1}{2} f'^2=0, f'\le 0 \}.
\end{equation*}
By the calculation in Theorem \ref{region_r1_thm}, the flow along $A$ is
inward-going.  Additionally, the flow along the first connected
component of $B$ is outward-going.  Finally, we put
$S(f,f',x)=\frac{1}{3}f^3 - \frac{1}{2} f'^2$ and observe that $\nabla
S$ is an inward pointing normal vector field to $B$.  We compute
\begin{eqnarray*}
\nabla S \cdot V &=& \begin{pmatrix} f^2 \\ - f' \\ 0 \end{pmatrix}^T 
\begin{pmatrix} f' \\ f^2 - \phi(x) \\ 1 \end{pmatrix} \\
&=& f' \phi(x) \le 0,
\end{eqnarray*}
so the flow along this component of $B$ is outward-going.  As a
result, we can apply the Antifunnel theorem, noting that $A$ is
connected, while $B$ is not.  Therefore, there exists a solution to
\eqref{ode1} that remains in $R_2$.  Note that there is a lower bound
on the $x$-coordinate of this solution, since the $x$-component of
$V(f,f',x)$ is equal to 1, and the Region $R_2$ lies within the
half-space $x>0$.  So this solution must enter $R_2$ through $A$, and
then never intersect $B$.  Additionally, notice that such a solution
will have $f' \le 0$ and $f \ge 0$, so it must be bounded.
\end{proof}
\end{lem}

\begin{lem}
\label{unbounded_funnel_lem}
Suppose $\phi(x)\ge 0$ and $\phi'(x) < 0$ for all $x \ge 0$.  The
complement of the set $A=R_1 \cup R_2$ consists of solutions which are
unbounded, and blow up in finite $x$. 
\begin{proof}
Let the complement of the set $A$ be called $C$, namely
$C=\{(f,f',x)|x>0\} - A$.  Now the calculations in Lemmas
\ref{region_r1_lem} and \ref{region_r2_lem} show that $C$ is a funnel,
in that the flow through the entire boundary of $C$ is inward.  If
$\phi$ does not tend to zero, then the argument in the proof of
Theorem \ref{limzero_iff_bounded_lem} completes the proof, as there is
a tubular neighborhood about $\{f=f'=0\}$ with strictly positive
radius in which solutions in $C$ cannot remain.  So without loss of
generality, we assume $\phi \to 0$.

\begin{figure}
\includegraphics[height=3in]{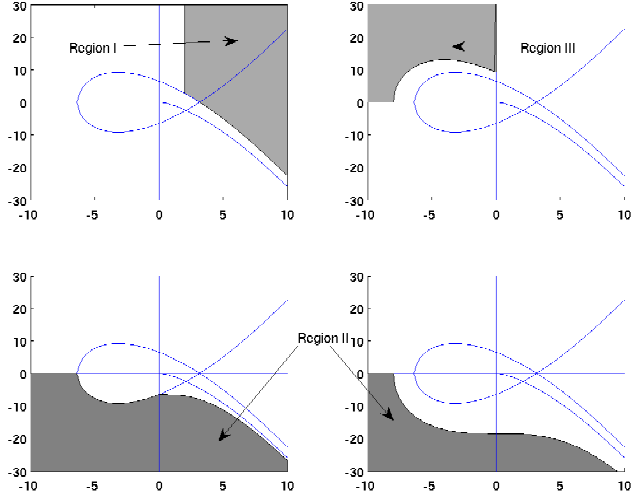}
\caption{The Regions $I$, $II$, and $III$ of Lemma \ref{unbounded_funnel_lem}}
\label{unboundedfunnel_fig}
\end{figure}

Define the Region $I$ by
\begin{equation*}
I=\left\{(f,f',x)| f > \sqrt{\phi(x)} \text{ and } \left(f'>0 \text{ or
  } H(f,f',x)>0\right) \right\}.
\end{equation*}
There are two bounding faces of Region $I$, along which the flow is
inward.  The first is $S_1 = f-\sqrt{\phi(x)} = 0$, along which  
\begin{eqnarray*}
\nabla S_1 \cdot V(f,f',x) &=& \begin{pmatrix}1\\ 0 \\ -
  \frac{\phi'(x)}{2 \sqrt{\phi(x)}} \end{pmatrix}^T
  \begin{pmatrix}f'\\ f^2 - \phi(x) \\ 1\end{pmatrix}\\
&=& f' - \frac{\phi'(x)}{2 \sqrt{\phi(x)}} > 0.
\end{eqnarray*}
The second was computed already in the proof of Theorem
\ref{region_r1_thm}.  Notice that $f''=f^2 - \phi(x) > 0 $ in Region
$I$, so $f(x)$ is concave-up, so solutions which enter Region $I$ are
unbounded.  Using similar reasoning to that of Theorem
\ref{limzero_iff_bounded_lem}, such solutions blow up in finite $x$.

Now suppose we have a point $(a,a',x_0) \in C$ with $a' < 0$.  We claim
that for some $x_1 > x_0$, the integral curve through this point will
cross the $f'=0$ plane.  To see this, construct Region $II$ by
\begin{equation*}
II=\left\{(f,f',x)|\frac{1}{3}f^3 - \frac{1}{2}f'^2 - \frac{1}{3}a^3 +
\frac{1}{2} a'^2 \le 0 \text{ and } f' \le 0\right\} \cap C.
\end{equation*}
Note that
\begin{equation*}
\nabla\left(\frac{1}{3}f^3 - \frac{1}{2} f'^2\right)\cdot V(f,f',x) =
f'\phi(x) \le 0,
\end{equation*}
so the flow is inward along Region $II$ except along $f'=0$ (along
which it is outward).  Also note that Region $II$ excludes a tubular
neighborhood of the line $f=f'=0$ with strictly positive radius.  As a
result of this, the integral curve though $(a,a',x_0)$ proceeds at
least as far as to allow $f < - \sqrt{\phi(x)}$, at which point, a
finite amount of distance in $x$ takes it to $f'=0$.

So at that point, the integral curve has entered Region $III$, say at
$x=x_1$, where
\begin{equation*}
III=\{(f,f',x)| H(f,f',x_1) \le 0 \text{ and } f \le 0 \text{ and } f'
\ge 0\}.
\end{equation*}
The flow is evidently inward along $f'=0$ and the curved portion by
previous calculations, and outward along $f=0$.  Again, note that the
line $f=f'=0$ is excluded from Region $III$ by a tubular neighborhood
of strictly positive radius, so there is an $x_2 > x_1$ where the
integral curve exits Region $III$ through $f=0$.

Now, consider a point $(0,c',x_2)$ along this integral curve with
$c'>0$.  In this case, the flow moves such a point rightward.  On the
other hand, the left boundary of Region $I$ moves leftward,
approaching $f=0$.  So there must be an $x_3 > x_2$ such that the
integral curve through $(0,c',x_2)$ enters the Region $I$.  Collecting
our findings, we see that every point in $C$ has an integral curve
which passes to Region $I$, and therefore corresponds to a solution
which is unbounded, and blows up for some finite $x$.  
\end{proof}
\end{lem}

\begin{thm}
\label{properties_of_Z}
Suppose $\phi(x)\ge 0$ and $\phi'(x) < 0$ for all $x \ge 0$.  The set $Z$
of initial conditions to \eqref{ode1} that lead to bounded solutions
\begin{enumerate}
\item lies within $A=R_1 \cup R_2$ and is
\item nonempty,
\item closed,
\item unbounded,
\item connected, and
\item simply connected.
\item Additionally, the portion of $Z$ corresponding to solutions that enter
the interior of $R_2$ is a 1-dimensional submanifold of $\{(f,f',x)|x=0\}$.
\end{enumerate}
\begin{proof}
\begin{enumerate}
\item From Lemma \ref{unbounded_funnel_lem}, all bounded solutions
  must lie in $A$. 
\item That there exist bounded solutions in $A$ is the content of
  Lemmas \ref{region_r1_lem} and \ref{region_r2_lem}.

\item Now, put $A_0 = A \cap \{(f,f',0)\}$ and $B_0 = \partial A -
  A_0$.  Observe that from the proofs of the previous theorems, the
  flow of $V$ along $A_0$ is inward, and the flow along $B_0$ is
  outward.  Since the last component of $V$ does not vanish, the flow
  of $V$ causes each point of $B_0$ to lie on an integral curve
  starting on $A_0$.  This establishes a homeomorphism $\Omega$ from
  $B_0$ into a subset of $A_0$.  In particular, $\Omega$ is an open
  map.  Now every solution passing through $B_0$ is of course
  unbounded, so $Z= A_0 - \Omega(B_0)$ is evidently closed (it is the
  complement of an open set).

\item $B_0$ clearly has the topology of $\mathbb{R} \times
  [0,\infty)$, so $\pi_1(B_0) = 0$.  Hence, $\pi_1(\Omega(B_0))=0$
  also, but notice that $\Omega(B_0)$ contains $\partial A_0$.  Suppose
  $Z$ were a bounded set.  Then it is contained in some disk $D$.  But
  $\partial D$ is homotopic to a loop in $A_0 - Z$, which either lies
  in $\text{int}(A_0 - Z)$ (in which case the homotopy need not move
  it) or in $\partial A_0$.  But this means that the loop encloses all
  of $Z$, and so cannot be contractible in $\Omega(B_0)$, which
  contradicts the triviality of $\pi_1(\Omega(B_0))$.   Hence $Z$ is
  unbounded.

\item We first show that the portion of $Z$ lying in the region $R_2$
satisfies the horizontal line test. First, note that a solution
starting in $Z\cap R_2$ cannot exit $R_2$.  For one, it cannot enter
$R_1$, since $R_1$ is an antifunnel.  Secondly, it cannot exit into
$\mathbb{R}^3 - (R_1 \cup R_2)$ since solutions there are all
nonglobal.  Suppose that $f_1(0) \ge f_2(0) \ge 0$ and
$f_1'(0)=f_2'(0)$ with $(f_1(0),f_1'(0))$ and $(f_2(0),f_2'(0))$ both
in $Z \cap R_2$.  But then
\begin{eqnarray*}
\frac{d}{dx}(f_1'(x)-f_2'(x))&=&f_1''(x)-f_2''(x)\\
&=&f_1^2(x) - f_2^2(x) \ge 0,\\
\end{eqnarray*}
with equality only if $f_1(0)=f_2(0)$.  Hence,
$\frac{d}{dx}(f_1(x)-f_2(x)) \ge 0$ for $x>0$, again with equality
only if $f_1(0)=f_2(0)$.  Now all solutions which remain in $R_2$ are
monotonic decreasing and bounded from below, so they must have
limits.  On the other hand, the only possible limit is
$(0,\lim_{x\to\infty}\sqrt{\phi(x)})$, so therefore all bounded
solutions in $R_2$ must have a common limit.  Therefore, we must
have that $f_1(0)=f_2(0)$.  Now this means that the portion of $Z$ in
the region $R_2$ can be realized as the graph of a function from the $f'$
coordinate to the $f$ coordinate.  Therefore, if $Z$ were not
connected, at least one component of $Z$ would be a bounded subset, which is
a contradiction.

\item Finally, if $Z$ were not simply connected, the Jordan curve
  theorem gives that there are two (or more) path components to
  $\Omega(B_0)=A_0 - Z$, which contradicts the continuity of $\Omega$.

\item By the connectedness of $Z$ and the horizontal line test in
  $R_2$, the function from the $f'$ coordinate to the $f$ coordinate
  whose graph is $Z\cap\text{int }R_2$ must be continuous.
  Additionally, by the connectedness of $Z$ and the uniquenss of
  solutions to ODE, this implies that the rest of $Z$ whose solutions
  enter the interior of $R_2$ is also a 1-manifold.
\end{enumerate}
\end{proof}
\end{thm}

\begin{df}
It is convenient to define, in addition to the initial condition set
$Z$, other sets $Z_{x_0} \subset \{(f,f',x)|x=x_0\}$ such that any
integral curve passing through a point in $Z_{x_0}$ exists for all
$x>0$.  Similarly, one can define $Z_{x_0}'$.
\end{df}

\begin{rem}
\label{approx_Z_rem}
If $\phi\to 0$ as $x\to\infty$, we conjecture that $Z$ acquires the
structure of a 1-manifold with boundary.  The series solution
\eqref{series_soln_coef} is not valid at such a boundary of $Z$, since
such a solution must remain in $R_1$ and therefore decays quicker than
the leading coefficient of \eqref{series_soln_coef}.  Indeed, by
analogy with the case where $\phi \equiv 0$, the leading term $f_0$ of
the series solution would vanish, and the solution is then asymptotic
to $-\int_x^\infty \int_t^\infty \phi(s) ds\, dt$.

All solutions in the form of the series solution
\eqref{series_soln_coef} enter $R_2$, so a result of this theorem is
that one of the two parameters $d$ or $K$ in the series solution is
superfluous.  Since $d$ parametrizes solutions when $\phi\equiv 0$, we
conventionally take $K=0$.  Using this, \eqref{handy_asymp_exp}
indicates that a good approximation (as $x_0 \to \infty$, locally near
$f=f'=0$) to the set $Z_{x_0}$ is the set
\begin{equation*}
\{H(f,f')=0\}=\{(f,f')|\frac{1}{3}f^3=\frac{1}{2}f'^2\}.
\end{equation*}
\end{rem}

\begin{rem}
If $\phi \to P >0$ as $x\to\infty$, then it is not true that $Z$ is a
1-manifold (with boundary).  Indeed, $Z$ has the structure of a
1-manifold attached to the teardrop-shaped set $M$ from Lemma
\ref{bounded_in_funnel_lem}.
\end{rem}

\section{Geometric properties of the initial condition set $Z$}
\label{geom_props_Z}

\begin{lem}
\label{Z_intersects_f_axis}
Suppose $\phi(x)>0$, $\phi'(x)<0$ for all $x>0$ and $\phi \to 0$ as
$x \to \infty$.  Then the set $Z$ intersects $\{(f,f',x)|f'=0\}$.
\begin{proof}
First, observe that $Z$ intersects the boundary of $R_1$ in $x=0$, since we
have by Lemmas \ref{region_r1_lem} and \ref{region_r2_lem} solutions
entirely within $R_1$ and its complement.  Using the fact that $Z$ is
connected and the Jordan curve theorem, $Z$ must intersect the
boundary of $R_1$ in the plane $x=0$.  This reasoning also applies for
each $Z_{x_0}$ with $x_0\ge 0$, so that we can find points in the
intersections $Z_{x_0} \cap \partial R_1$ for each $x_0 \ge 0$.  Also
note that for the backwards flow associated to our equation (ie. the
flow of $-V$), solutions which enter $R_1$ must exit through the plane
$x=0$.  Hence there exists a sequence of points $\{F_n\} \subset Z$
with $F_n=(f_n,f_n',0)$ such that the integral curve through $F_n$
passes through $G_n=(g_n,g_n',n) \in Z_n \cap \partial R_1$ for each
integer $n \ge 0$.

Discern three cases:
\begin{enumerate}
\item If any $F_n$ are in Quadrants I or II, then since $Z$ is
  connected, it must intersect $\{f'=0\}$.
\item If any $F_n$ are in Quadrant III, observe that the flow across
  the surface $S=\left \{(f,f',x)|\frac{1}{3}f^3=\frac{1}{2}f'^2, f' \le
  0\right \}$ is right-to-left.  Thus the integral curve must cross
  into Quadrant II on its way to $G_n$.  Therefore, the set $Z$
  cannot intersect the surface $S$, and so it must intersect
  $\{f'=0\}$.
\item Assume all the $F_n$ lie in Quadrant IV.  Observe that $\{F_n\}$
  is a closed subset of $R_1 \cap \{x=0\}$, which is compact.
  Hence some subsequence of $\{F_n\}$ must have a limit, say $F$.
  Since $Z$ is closed, $F \in Z$.  But in the portion of $R_1$ lying
  in the $x=0$ plane and in Quadrant IV, we have that 
\begin{equation*}
\frac{d}{dx}f'=f^2 - \phi < 0
\end{equation*}
and
\begin{equation*}
\frac{d}{dx}f=f' < 0.
\end{equation*}
Hence $f_n' \ge g_n'$.  But since $\phi \to 0$, $g_n' \to 0$, so $F$
lies on $\{f'=0\}$.
\end{enumerate}
\end{proof}
\end{lem}

\begin{lem}
\label{Z_intersects_fp_axis_weak}
Under the same hypotheses as Lemma \ref{Z_intersects_f_axis}, $Z$ also
intersects the half plane $\{f=0,f'>0\}$.
\begin{proof}
Using Lemma \ref{Z_intersects_f_axis}, we form a sequence $\{F_n\}
\subset Z$ such that the integral curve through $F_n$ passes through
$\{f'=0,f \ge 0, x=n\}$ for each integer $n$.  (This can be done
without loss of generality, because if any integral curves pass
through $\{f'=0,f < 0\}$, then the proof is complete by connectedness
of $Z$.)  Note that this sequence is entirely contained within $R_1$
by Lemma \ref{unbounded_funnel_lem}.

Discern three cases:
\begin{enumerate}
\item There exists an $F_n$ in either of Quadrants II or III.  The
  result follows by the connectedness of $Z$.
\item There exists $F_n$ in Quadrant IV.  This cannot occur unless the
  integral curve through $F_n$ passes through Quadrant III since the
  flow along $\{f'=0\}$ points inward into the portion of Quadrant IV
  inside $R_1$.
\item Otherwise, we assume $\{F_n\}$ is entirely contained within
  Quadrant I.  In this case, note that 
\begin{equation*}
\frac{d}{dx} f' = f^2 - \phi < 0.
\end{equation*}
Hence the $f'$-coordinate of the integral curve through each $F_n$ is
positive on the interior of Quadrant I.  Hence 
\begin{equation*}
\frac{d}{dx} f = f' > 0,
\end{equation*}
so $f_n \le g_n$.  But $g_n \to 0$ since $\phi \to 0$, so any limit
point of $\{F_n\}$ will have $f$-coordinate equal to zero.  By the
compactness of $R_1 \cap \{x=0\}$ and the closedness of $Z$, this
implies that $Z$ intersects $\{f=0,f'>0\}$.
\end{enumerate}
\end{proof}
\end{lem}

\begin{lem}
\label{Z_intersects_fp_axis}
Suppose $\phi(x) >0$ for all $x>0$, $\phi \to 0$ as $x \to \infty$,
and that there exists an $x_0 \ge 0$ such that for all $x>x_0$,
$\phi'(x)<0$.  Then the set $Z$ intersects $\{f=0, f'>0\}$.
\begin{proof}
We follow the pattern of proving the existence of an intersection for
an open interval in $x$ containing $x_0$, and then constructing an
{\it a priori} estimate for the $f'$-coordinate of this intersection.  

Apply Lemma \ref{Z_intersects_fp_axis_weak} to $x_0$, we have that
$Z_{x_0}$ intersects $\{f=0,f'>0\}$.  Let $(0,f_0',x_0)$ lie in this
intersection.  Note that 
\begin{equation*}
\frac{d}{dx} f = f'>0
\end{equation*}
and
\begin{equation*}
\frac{d}{dx}f' = f^2 - \phi = -\phi < 0
\end{equation*}
when evaluated there.  As a result, the integral curve passing through
$(0,f_0',x_0)$ must pass through Quadrant II first, say for $x \in
(x_1,x_0)$.  Then evidently, $Z_{x_1}$ must intersect $\{f=0,f'>0\}$.

Now since $\phi(x)>0$ between $x_1$ and $x_0$, and $[x_1,x_0]$ is
compact, there is an open set in $\mathbb{R}^3$ containing the
intersection of each $Z_x$ with $\{f=0,f'>0\}$ for each $x\in
[x_1,x_0]$, such that in this open set $\frac{d}{dx} f' \le K < 0$.  As a
result, $f_1' \ge f_0'$.  Hence the $f'$-coordinate of the
intersection point of $Z_x$ with $\{f=0,f'>0\}$ is decreasing with
increasing $x$.  (Since we have $f^2 - \phi > -\phi$, it is decreasing
at a rate no faster than $\phi$.  This implies that this intersection point
has $f'$-coordinate no larger than $\int_0^{x_0} \phi(x) dx + f_0'$ at
$x=0$.)  Now since solutions through $Z_{x_1}$ exist for all $x>0$ by
definition, this suffices to show that $Z$ intersects $\{f=0,f'>0\}$.
\end{proof}
\end{lem}

\begin{rem}
\label{delicateness_rem}
The line of reasoning used in the third case of each of Lemmas
\ref{Z_intersects_f_axis} and \ref{Z_intersects_fp_axis_weak} (and
also in \ref{Z_intersects_fp_axis}) fails if we try to continue $Z$ much
farther.  This is due to the nonmonotonicity of $df'/dx$ in Quadrants
II and III.  More delicate control of $\phi$ must be exercised to say
more.
\end{rem}

\begin{calc}
Towards the end of the more delicate results mentioned in Remark
\ref{delicateness_rem}, it is useful to know the maximum speed along
integral curves on points in the region $R_1$ in the $f$- and
$f'$-directions.  By this we mean to compute for fixed $x$ the maximum
values of 
\begin{equation}
\label{directional_speeds}
\begin{cases}
|f'| \text{ for the }f\text{-direction}\\
|f^2-\phi(x)| \text{ for the }f'\text{-direction}\\
\end{cases}
\end{equation}
in $R_1$.  The first is easy to maximize: we simply look for the
maximum value of $f'$ in $R_1$, which is a maximum of 
\begin{equation*}
f'= \sqrt{\frac{2}{3} f^3 - 2f\phi(x) + \frac{4}{3} \phi^{3/2}(x)},
\end{equation*}
for $-2\sqrt{\phi(x)} \le f \le \sqrt{\phi(x)}$.  This occurs at
$f=-\sqrt{\phi(x)}$, and has the value of $\sqrt{8/3}\phi^{3/4}$.
For the second part of \eqref{directional_speeds}, it is easy to see
that the maximum is $3 \phi(x)$. In summary,
\begin{equation}
\label{max_speed_calc}
\begin{cases}
|f'| \le \sqrt{\frac{8}{3}}\phi^{3/4}(x) \text{ for the }f\text{-direction}\\
|f^2-\phi(x)| \le 3\phi(x) \text{ for the }f'\text{-direction}\\
\end{cases}
\end{equation}
on $R_1$.
\end{calc}

Using this calculation, we can impose a stronger bound on the decay of
$\phi(x)$, and constrain the set $Z$ further.

\begin{lem}
\label{constrain_Z_lem_weak}
Suppose $\phi(x)>0$, $\phi'(x) < -D
\frac{4\sqrt{2}}{k\sqrt{3}}\phi^{5/4}(x)$ for all $x>0$ for some
$0<k<1$ and $D>1$.  Then the set $Z$ is contained within $\{f \ge
-k\sqrt{\phi(0)}\}$ and intersects each vertical and horizontal line
in $\{f \ge 0\}$ exactly once, and intersects $\{f'=0\}$ only once.
\begin{proof}
That $Z$ intersects $\{f=0,f' \ge 0\}$ and $\{f'=0,f\ge 0\}$ at all
follows from Lemmas \ref{Z_intersects_f_axis} and
\ref{Z_intersects_fp_axis}.  Now consider the region $A \subset R_1$
shown in Figure \ref{constrain_Z_fig} and defined by 
\begin{equation*}
A =R_1 \cap \left(\{f'\ge 0, f \le k \sqrt{\phi(x)}\} \cup \{f' \le
0, 2 f^3 \le 3 f'^2\} \right).
\end{equation*}

\begin{figure}
\includegraphics[height=3in]{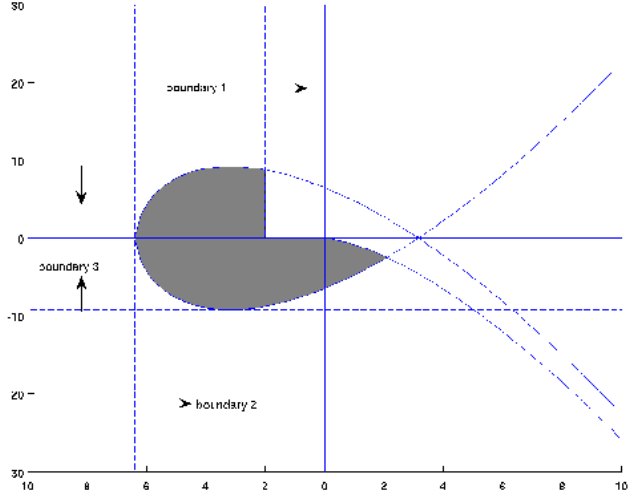}
\caption{The region $A$ of Lemma \ref{constrain_Z_lem_weak}}
\label{constrain_Z_fig}
\end{figure}

The boundary segments strictly to the right of the boundary labelled 1 in
Figure \ref{constrain_Z_fig} are evidently inflow, so long as $\phi >
0$.  The boundary labelled as 1 in the figure moves with speed
\begin{equation*}
\frac{d}{dx} ( - k \sqrt{\phi(x)}) = \frac{-k}{2\sqrt{\phi(x)}}
\phi'(x) > D \frac{2\sqrt{2}}{\sqrt{3}} \phi^{3/4}(x)
\end{equation*}
which is greater than maximum speed in the $f$-direction given in
\eqref{max_speed_calc}.  This implies that the boundary moves faster
than any solution inside $R_1$.  Hence it is an inflow portion of the
boundary.  On the other hand, the curved segment of the boundary to
the left has been shown to be outflow, in Lemma \ref{region_r1_lem}.

We observe that the boundary marked 2 in Figure
\ref{constrain_Z_fig} moves with speed
\begin{equation*}
\frac{d}{dx} ( - 2 k \sqrt{\phi(x)} ) = \frac{ -k}{\sqrt{\phi(x)}} \phi'(x),
\end{equation*}
which is strictly faster than the boundary marked 1 in Figure
\ref{constrain_Z_fig}, and the boundary marked 3 in Figure
\ref{constrain_Z_fig} moves with speed
\begin{equation*}
\frac{d}{dx} \left( \pm \sqrt{\frac{8}{3}} \phi^{3/4}(x) \right ) =
\pm \frac{\sqrt{3}}{\sqrt{2} \phi^{1/4}(x)} \phi'(x),
\end{equation*}
noting that $f'(-\sqrt{\phi(x)})=\pm\sqrt{8/3}\phi^{3/4}(x)$ is the
value of the maximum $f'$-coordinate of $R_1$ at a given $x$ value.
This last speed is greater than the maximum speed in the $f'$-direction given by
\eqref{max_speed_calc} since $\phi'(x) < - D 
\sqrt{6}\phi^{5/4}(x)$.  (Notice that
$\sqrt{6} < \frac{4\sqrt{2}}{k\sqrt{3}}$, since $0<k<1$.)  

Since $D>1$, this means that both the boundaries marked 2 and 3 in
Figure \ref{constrain_Z_fig} overtake any solution constrained to be
within $R_1$.  As a result, every solution within the region $A$ must
leave it within finite $x$.  But the only way to leave $A$ causes a
solution to enter $\mathbb{R}^3 - (R_1 \cup R_2)$, so every solution
which contains a point in $A$ cannot exist for all $x>0$ by Lemma
\ref{unbounded_funnel_lem}.  Therefore, $Z$ is contained within $(R_1
\cup R_2) - A$.

Now consider the region $B$ which is defined by
\begin{equation*}
B =R_1 \cap \left(\{f'\ge 0, f \le 0\} \cup \{f' \le
0, 2 f^3 \le 3 f'^2\} \right),
\end{equation*}
which is simply the region $A$, with $k$ taken to be zero.  The
portion of the boundary of $B$ lying in the $\{f=0\}$ plane is inflow.
We can therefore apply the reasoning of the vertical line test:
Suppose $(f_1,f_1',0),(f_2,f_2',0) \in Z$ with $f_1=f_2 > 0$ and $f_1'
\ge f_2'$.  Then we have both (at $x=0$)
\begin{equation*}
\frac{d}{dx}(f_1'-f_2')=f_1^2-f_2^2 =0 
\end{equation*}
and
\begin{equation*}
\frac{d}{dx}(f_1-f_2)=f_1'-f_2' \ge 0,
\end{equation*}
which gives that $f_1^2-f_2^2 \ge 0$ for some open interval about
$x=0$.  Then, $\frac{d}{dx}(f_1'-f_2')\ge 0$, which implies that
in fact $\frac{d}{dx}(f_1-f_2) \ge 0$.  However, since all bounded
solutions tend to the common limit of zero, we have that this implies
$f_1'=f_2'$ at $x=0$.  (Note that since each solution starts in $Z \cap
(R_1 - B)$, we have that neither solution can become negative, since
that would involve entering $B \subset A$ or leaving $R_1 \cup R_2$.)
This implies that there is a unique intersection of $Z$ with each vertical
line.  The same reasoning applies in the case of the horizontal line
test, as in Theorem \ref{properties_of_Z}.
\end{proof}
\end{lem}

\begin{lem}
\label{constrain_Z_lem}
Suppose $\phi(x)>0$ for all $x>0$, and that $\phi'(x) < -D 
\frac{4\sqrt{2}}{k\sqrt{3}}\phi^{5/4}(x)$ for all $x>x_0\ge 0$ for some $0<k<1$ and $D>1$.
Additionally, suppose that for all $x \in [0,x_0]$, 
\begin{equation}
\label{constrain_Z_eqn}
x_0 - x < \frac{\sqrt{\phi(x)}-k\sqrt{\phi(x_0)}}{
\sqrt{\frac{8}{3}} P^{3/4}},
\end{equation}
where $P=\max_{x \in [0,x_0]} \phi(x)$. Then the set $Z$ is contained
within $\{f \ge -\sqrt{\phi(0)}\}$ and intersects each vertical and
horizontal line in $\{f \ge 0\}$ exactly once, and intersects
$\{f'=0\}$ only once.
\begin{proof}
The set $Z_{x_0}$ is constrained to lie within the set $\{f' \ge
-k\sqrt{\phi(x_0)}\}$, by Lemma \ref{constrain_Z_lem_weak} (replacing
$x_0$ by zero).  Now using the $f$-direction part of
\eqref{max_speed_calc}, the smallest $f$-value attained in $Z_x$ is
\begin{equation*}
\int_{x_0}^x \sqrt{\frac{8}{3}}\phi^{3/4}(x) dx -
k \sqrt{\phi(x_0)}.
\end{equation*}
If $x<x_0$, we have
\begin{eqnarray*}
\int_{x_0}^x \sqrt{\frac{8}{3}}\phi^{3/4}(x) dx -
k \sqrt{\phi(x_0)} &\ge& \sqrt{\frac{8}{3}} P^{3/4} (x-x_0) - k\sqrt{\phi(x_0)}\\
&>&-\sqrt{\phi(x)},\\
\end{eqnarray*}
by \eqref{constrain_Z_eqn}.  As a result, $Z_x \subset \{f \ge
-\sqrt{\phi(x)}\}$ for each $x<x_0$.  This additionally means that in
the backwards flow, the entire portion of $Z_x$ contained in $\{f \le
0\}$ is moving away from the plane $\{f'=0\}$, which completes the
proof.
\end{proof}
\end{lem}

\begin{rem}
The condition that $\phi'(x) < C \phi^{5/4}(x)$
implies 
\begin{eqnarray*}
\phi^{-5/4} \phi'(x) &<& C\\
-\frac{1}{4} \phi^{-1/4}(x) &<& Cx + C'\\
\phi(x) &<&\frac{C'''}{\left(C''-x\right)^4}
\end{eqnarray*}
for some $C''$ and $C'''$.  Notice that this condition is satisfied
when the series solution converges by Theorem \ref{series_conv_thm}.
\end{rem}

\section{Solutions on the entire real line}
\label{long_sec}
We now combine the results for \eqref{ode1} and \eqref{ode1_backwards}
to discuss properties of the solutions to \eqref{long_ode1}.
When $\phi(x)$ is monotonically decreasing, we have by Lemma
\ref{unbounded_funnel_lem} that the initial condition set stays within
$R_1 \cup R_2$.  In particular, $Z \subset \{f' \le \sqrt{\frac{8}{3}}
\phi^{3/4}(0)\}$.  If we relax the restriction of monotonicity, we
obtain a similar result.

\begin{lem}
\label{bound_on_fp}
If $f=f(x)$ is a bounded solution to the initial value problem
\eqref{ode1} with $\phi \in C^\infty\cap
  L^\infty(\mathbb{R})$ then $f'(0) < \sqrt{8/3}
    \|\phi\|_\infty^{3/4}$.
\begin{proof}
Since $f$ is a solution to \eqref{ode1}, then it must satisfy
\begin{equation*}
f''=f^2-\phi(x) \ge f^2 - \|\phi\|_\infty.
\end{equation*}
Now Lemma \ref{bounded_in_funnel_lem} shows that all bounded solutions
to $g''=g^2 - \|\phi\|_\infty$ lie in the closure of the set $M$ given
by
\begin{equation*}
M=\left \{(g,g')|\frac{1}{3} g^3 - \frac{1}{2} g'^2 - g
\|\phi\|_\infty + \frac{2}{3} \|\phi\|_\infty^{3/2}>0,g <
\sqrt{\|\phi\|_\infty} \right\}.
\end{equation*}

Since this set $M$ is bounded, we can find the maximum value of $f'$,
which is $f'_{max} = \sqrt{8/3}\|\phi\|_\infty^{3/4}.$
\end{proof}
\end{lem}

\begin{lem}
\label{bound_on_phi_integral_1}
Consider solutions to \eqref{long_ode1} on the real line, with $\phi
\in C_0^\infty \cap L^\infty(\mathbb{R})$.  If for some $-\infty < A <
B < \infty$,
\begin{equation*}
-\int_A^B \phi(x) dx > \sqrt{\frac{8}{3}} \left (
 (\sup_{x\in(-\infty,A]} |\phi(x)|)^{3/4} +
 (\sup_{x\in[B,\infty)} |\phi(x)|)^{3/4}  \right)
\end{equation*}
then no bounded solutions exist.
\begin{proof}
For a solution $f$, we have that $f''=f^2-\phi(x) \ge -\phi(x).$
Integrating both sides we have
\begin{equation*}
f'(B)-f'(A) \ge - \int_A^B \phi(x) dx.
\end{equation*}
By Lemma \ref{bound_on_fp}, bounded solutions on
\begin{itemize}
\item $x>B$ have $f'(B) < \sqrt{\frac{8}{3}} (\sup_{x\in(-\infty,A]}
  |\phi(x)|)^{3/4}$, and
\item on $x>A$, they have $f'(A) < (\sup_{x\in[B,\infty)}
  |\phi(x)|)^{3/4}$,
\end{itemize}
so a necessary condition for there to be a bounded solution is that
\begin{equation*}
-\int_A^B \phi(x) dx \le \sqrt{\frac{8}{3}} \left (
 (\sup_{x\in(-\infty,A]} |\phi(x)|)^{3/4} +
 (\sup_{x\in[B,\infty)} |\phi(x)|)^{3/4}  \right).
\end{equation*}
\end{proof}
\end{lem}

\begin{cor}
\label{phi_integral_bound}
A necessary condition for bounded solutions to \eqref{long_ode1} to exist if $\phi \in
C_0^\infty \cap L^\infty(\mathbb{R})$ is $\int_{-\infty}^\infty \phi(x) dx > 0$.
\begin{proof}
Suppose bounded solutions exist.  By the proof of Lemma
\ref{bound_on_phi_integral_1}, if we let
\begin{equation*}
g_n=-\int_{-n}^n \phi(x) dx,
\end{equation*}
and
\begin{equation*}
h_n=\sqrt{\frac{8}{3}}\left (
 (\sup_{x\in(-\infty,-n]} |\phi(x)|)^{3/4} +
 (\sup_{x\in[n,\infty)} |\phi(x)|)^{3/4}  \right),
\end{equation*}
then $g_n < h_n$ for each positive integer $n$.  But the continuity of
limits gives 
\begin{equation*}
- \int_{-\infty}^\infty \phi(x) dx = \lim_{n\to \infty} g_n < \lim_{n
  \to \infty} h_n = 0.
\end{equation*}
\end{proof}
\end{cor}

\begin{df}
A function $\phi \in C_0^\infty \cap L^\infty(\mathbb{R})$ will be
  called {\it M-shaped} if there exists an $x_0>0$ such that for all
  $|x|>x_0$, $\phi(x) > 0$ and
\begin{itemize}
\item $\phi$ is monotonic increasing for $x < -x_0$ and 
\item $\phi$ is monotonic decreasing for $x > x_0$.
\end{itemize}
\end{df}

\begin{thm}
\label{long_soln_exists}
Suppose $\phi$ is a positive M-shaped function, then solutions exist
to \eqref{long_ode1}. 
\begin{proof}
Observe that by Lemma \ref{Z_intersects_fp_axis}, we have that the set
$Z$ intersects $\{f=0,f'>0\}$.  Additionally, by Theorem
\ref{properties_of_Z}, we have that $Z$ also lies in $R_2$, which is
unbounded in Quadrant IV.  Likewise, the set $Z'$ (for
\eqref{ode1_backwards}) intersects $\{f=0,f'<0\}$, and becomes
unbounded in Quadrant I, so $Z \cap Z'$ must be nonempty,
and at least one point in this intersection is in the half-plane $\{x=0,f>0\}$.
\end{proof}
\end{thm}

\begin{thm}
\label{long_soln_unique}
Suppose $\phi$ is a positive M-shaped function which additionally
satisfies the decay constraints of Lemma \ref{constrain_Z_lem} for
$x>0$ and $x<0$ seperately, then a unique positive solution exists
to \eqref{long_ode1}.  (Note that for $x<0$, the inequalities and
signs in Lemma \ref{constrain_Z_lem} must be reversed, {\it mutatis
  mutandis}.)
\begin{proof}
By the Theorem \ref{long_soln_exists}, there exist solutions to
\eqref{long_ode1}, one of which comes from the intersection of $Z \cap
Z'$ in the half-plane $\{x=0,f>0\}$.  The vertical-line test in Lemma
\ref{constrain_Z_lem} allows one to conclude that the solution which
passes through that half-plane must continue directly to the region
$R_2$ of Lemma \ref{region_r2_lem}, without crossing the plane
$\{f=0\}$.  Thus this solution is strictly positive.

On the other hand, Lemma \ref{constrain_Z_lem} indicates that $Z$ may
lie only in Quadrants I, II, and IV, while the set $Z'$ must lie in
Quadrants I, III, and IV.  On the other hand, the vertical- and
horizontal-line tests ensure a unique intersection of $Z$ and $Z'$ in
Quadrants I and IV, so the solution is unique.
\end{proof}
\end{thm}

\begin{eg}
\label{example_1}
We examine the family $\phi_c(x)=c e^{-x^2/2}$, which is M-shaped when
$c\ge 0$.  Notice that when $c<0$, then the necessary condition of
Corollary \ref{phi_integral_bound} is not met, so solutions do not
exist for all $x \in \mathbb{R}$.  When $c=0$, then the trivial
solution $f=0$ is the only solution.  For $c>0$, we examine
$\phi'_c(x) = -xce^{-x^2/2}$.  Figure \ref{example_1} shows the sets
$Z$ and $Z'$ for the case when $c=0.05$.  In particular, one notes
that there appears to be a unique point of intersection.

\begin{figure}
\includegraphics[height=3in]{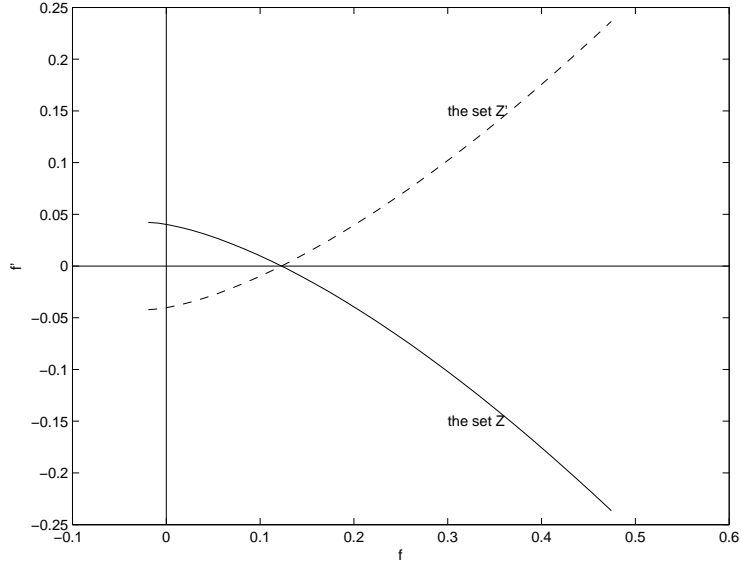}
\caption{The sets $Z$ and $Z'$ in Example \ref{example_1}}
\label{Z_example_1}
\end{figure}

We find the $x_0$ for which larger $x$ satisfy $\phi'(x) <
-4\sqrt{2}\phi^{5/4}(x)/(k\sqrt{3})$:
\begin{eqnarray*}
-xce^{-x^2/2} &<& -\frac{4\sqrt{2}}{k\sqrt{3}}c^{5/4}e^{-5x^2/8}\\
xe^{x^2/8}&>&\frac{4\sqrt{2}}{k\sqrt{3}}c^{1/4},
\end{eqnarray*}
which occurs if $x>\frac{4\sqrt{2}}{k\sqrt{3}}c^{1/4}$, so we may take
$x_0=\frac{4\sqrt{2}}{k\sqrt{3}}c^{1/4}$.  

By way of example, if we fix $x_0 = 4/3$, then $k=\sqrt{6}c^{1/4}$.
(We enforce $0<k<1$ by taking $c$ small.)  Now we must check to see if
\eqref{constrain_Z_eqn} holds.  In this case, we need to see if $c$
can be chosen so that $x_0-x=4/3-x$ is bounded above by
\begin{eqnarray*}
\frac{\sqrt{\phi(x)}-k\sqrt{\phi(x_0)}}{
\sqrt{\frac{8}{3}} P^{3/4}} &=&
\frac{\sqrt{c}e^{-x^2/4}-\sqrt{6} c^{3/4}
  e^{-16/36}}{\sqrt{8/3}c^{3/4}}\\
&=&\frac{e^{-x^2/4}-\sqrt{6}c^{1/4}e^{-16/36}}{\sqrt{8/3}c^{1/4}}\\
&\ge&\frac{e^{-16/36}-\sqrt{6}c^{1/4}e^{-16/36}}{\sqrt{8/3}c^{1/4}},
\end{eqnarray*} 
which can be made as large as one likes by taking $c$ sufficiently
small.  Noting that this last line is a constant in $x$ completes the
bound.  Therefore, there is a unique positive solution for
$0=f''-f^2+ce^{-x^2/2}$ with $c \in [0,\epsilon)$ for some
  $\epsilon>0$.

\end{eg}

\begin{rem}
\label{asymp_summary_rem}
Taken together, the results of Corollary \ref{phi_integral_bound} and
Theorems \ref{long_soln_exists} and \ref{long_soln_unique} for
M-shaped $\phi$ provide the following story about solutions to the
equation $0=f''-f^2+\phi$ on the real line:
\begin{itemize}
\item If the portion of $\phi$ where it is allowed to be negative is
  sufficiently negative, then no solutions exist,
\item If $\phi$ is positive, then a solution will exist.  There is no
  particular reason to believe that this solution will be strictly
  positive or unique.  
\item If the decay in the monotonic portions of $\phi$ is fast enough,
  there is exactly one solution, which is strictly positive.
\end{itemize}
\end{rem}

\section{Numerical examination}
\label{numer_sec}
\subsection{Computational framework}

Notice that the results of Remark \ref{asymp_summary_rem} are not
sharp: nothing is said if $\phi$ has a portion which is negative, but
still satisfies the necessary condition of Corollary
\ref{phi_integral_bound}.  Further, if $\phi$ is positive, but does
not satisfy the decay rate conditions, nothing is said about the
number of global solutions that exist.  Answers to these questions can
be obtained by combining the asymptotic information we have collected
about the sets $Z$ and $Z'$ with a numerical solver.  In particular,
we can obtain information about the number of global solutions to
\eqref{long_ode1} for any M-shaped $\phi$.

Suppose that $\phi$ is an M-shaped function, and that $x_0$ is such
that $\phi(x)$ is monotonic decreasing for all $x>x_0$ and is
monotonic increasing for all $x<-x_0$.  (If $\phi$ decreases fast
enough, we can choose $x_0$ so that the series solution converges on
the complement of $(-x_0,x_0)$ for sufficiently small initial
conditions.)  Then we have the sets $Z'_{-x_0}$ and $Z_{x_0}$ of
initial conditions to ensure existence of solutions on
$(-\infty,-x_0]$ and $[x_0,\infty)$ respectively.  Then any solution
to the boundary value problem
\begin{equation}
\label{ode1_bvp}
\begin{cases}
0=f''(x)-f^2(x)+\phi(x) \text{ for } -x_0<x<x_0\\
(f(-x_0),f'(-x_0))\in Z'_{-x_0}\\
(f(x_0),f'(x_0))\in Z_{x_0}\\
\end{cases}
\end{equation}
extends to a global solution of \eqref{long_ode1}.  So all one must do
is solve \eqref{ode1_bvp} numerically.  An easy way to do this is to
numerically extend the sets $Z'_{-x_0}$ and $Z_{x_0}$ to $Z'$ and $Z$
respectively (ie. extend them to the plane $x=0$) and compute $Z' \cap
Z$.  

\begin{figure}
\includegraphics[height=3in]{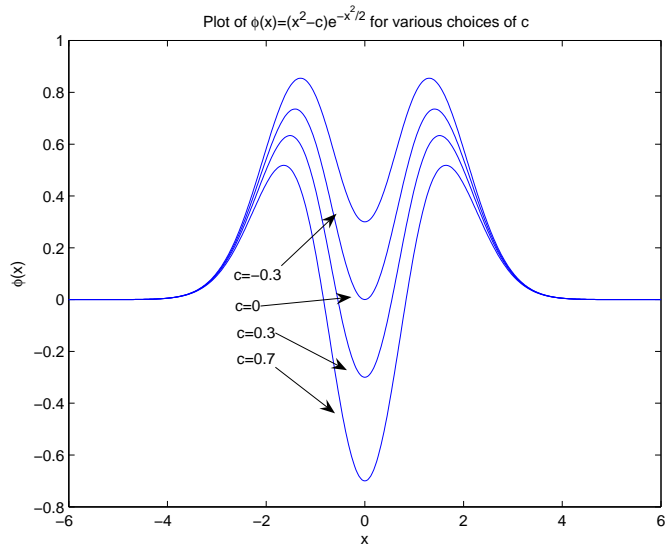}
\caption{The function $\phi(x;c)$ for various $c$ values}
\label{sample_phi_fig}
\end{figure}

In order to analyze \eqref{long_ode1} numerically, it is
necessary to make a choice of $\phi$.  Evidently, the numerical
results for that particular choice of $\phi$ cannot be expected to
apply in general.  However, a good choice of $\phi$ will suggest
features in the solutions that are common to a larger class of
$\phi$.  We shall use
\begin{equation}
\label{sample_phi}
\phi(x;c)=(x^2-c)e^{-x^2/2},
\end{equation}
where $c$ is taken to be a fixed parameter.  (See Figure
\ref{sample_phi_fig})  This choice of $\phi$ has the following
features which make for interesting behavior in solutions to
\eqref{long_ode1}:
\begin{itemize}
\item $\phi(x;c)>0$ for $c<0$.  In this case, there are solutions to
  \eqref{long_ode1}, by Theorem \ref{long_soln_exists}.  On the other
  hand, the decay rate conditions are not met over all of $\mathbb{R}$
  so the uniqueness result of Theorem \ref{long_soln_unique} does not
  apply.  Inded, the decay rate conditions are met only for
  sufficiently large $|x|$, but not for $|x|$ small.  
\item If $c>0$ is large enough, it should happen that no solutions to
  \eqref{long_ode1} exist, since the necessary condition of Corollary
  \ref{phi_integral_bound} is not met.  Indeed, the integral of $\phi$
  vanishes when $c=1$.
\end{itemize}

\subsection{Bifurcations in the global solutions}

\begin{figure}
\includegraphics[height=2in]{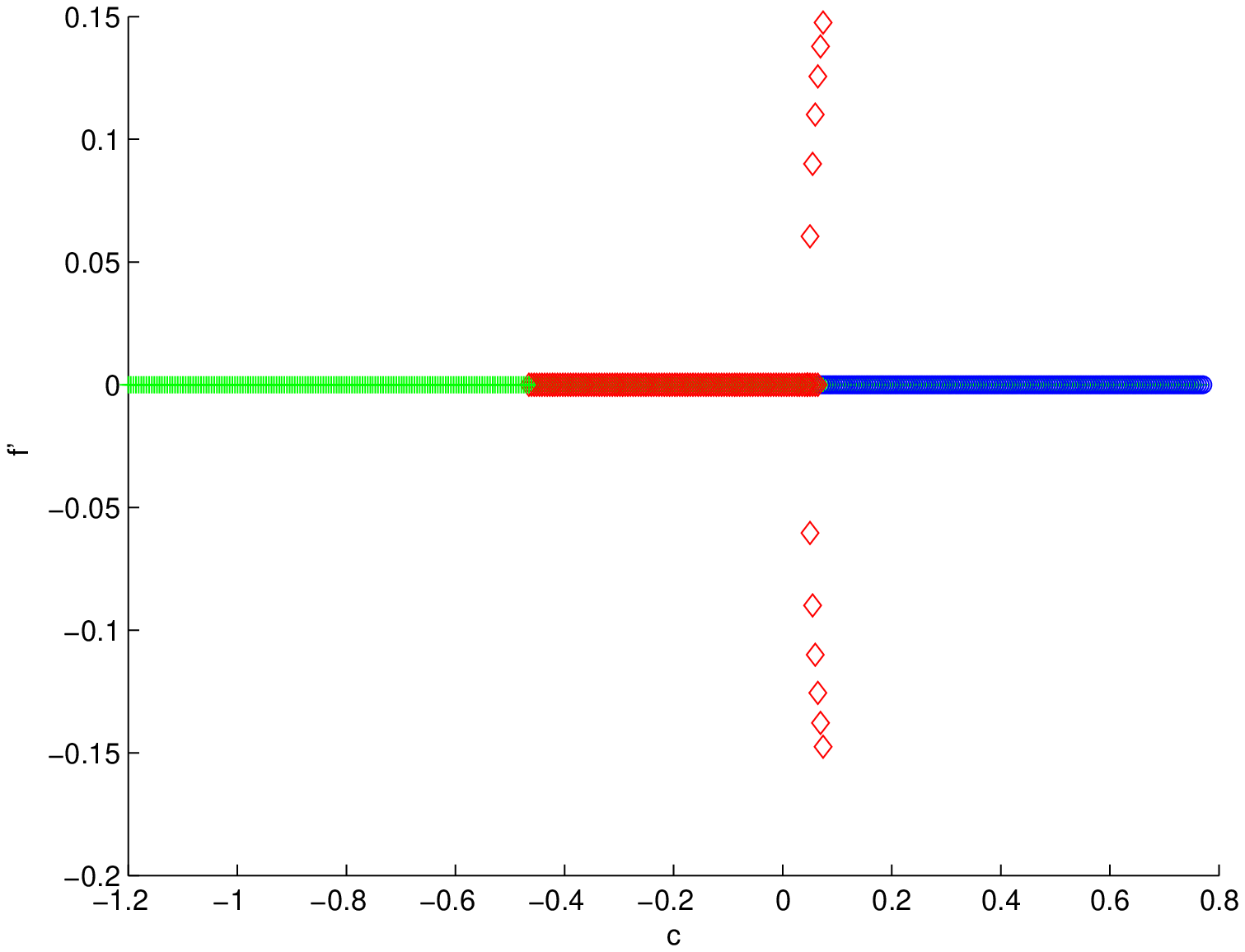}
\includegraphics[height=2in]{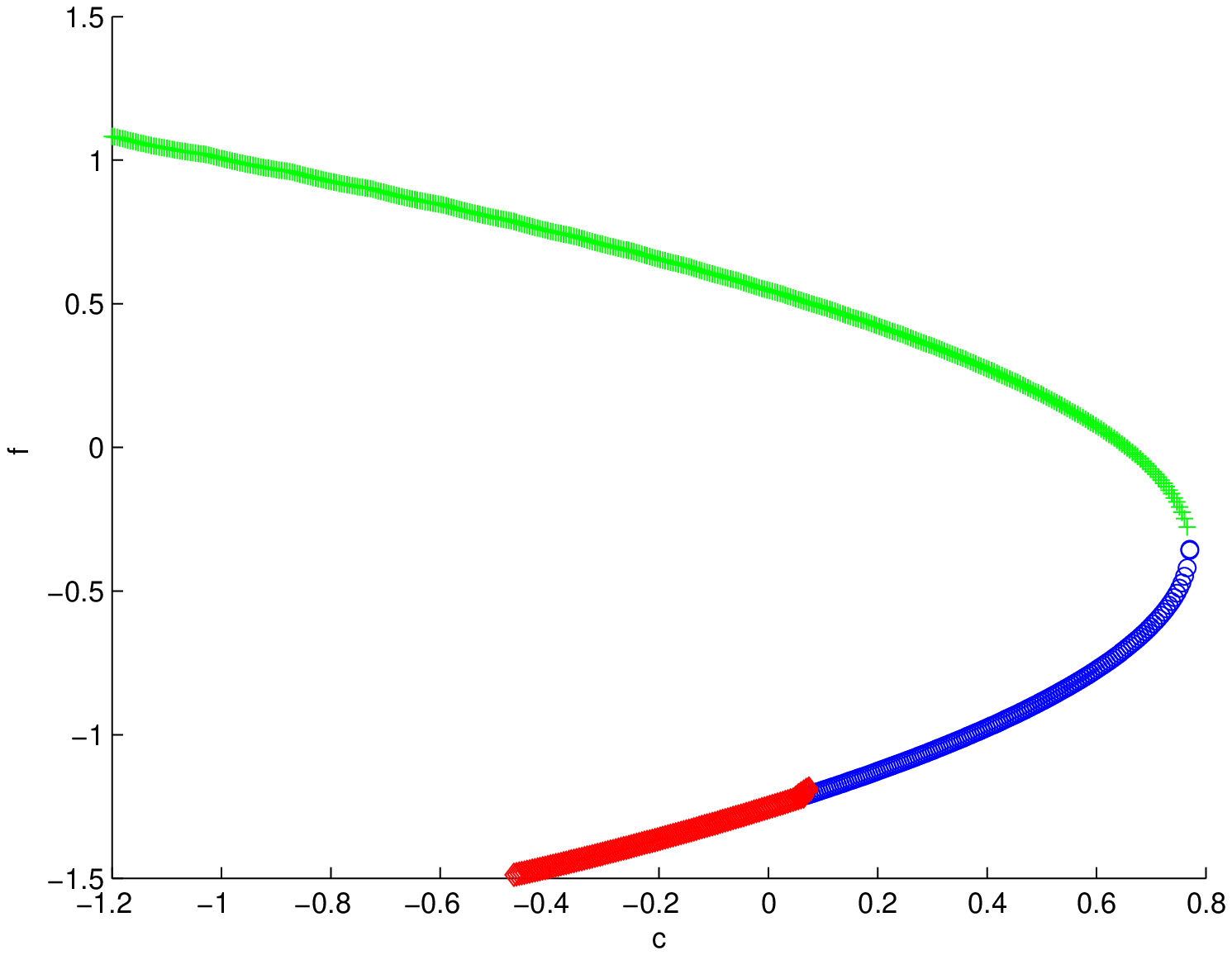}
\includegraphics[height=2in]{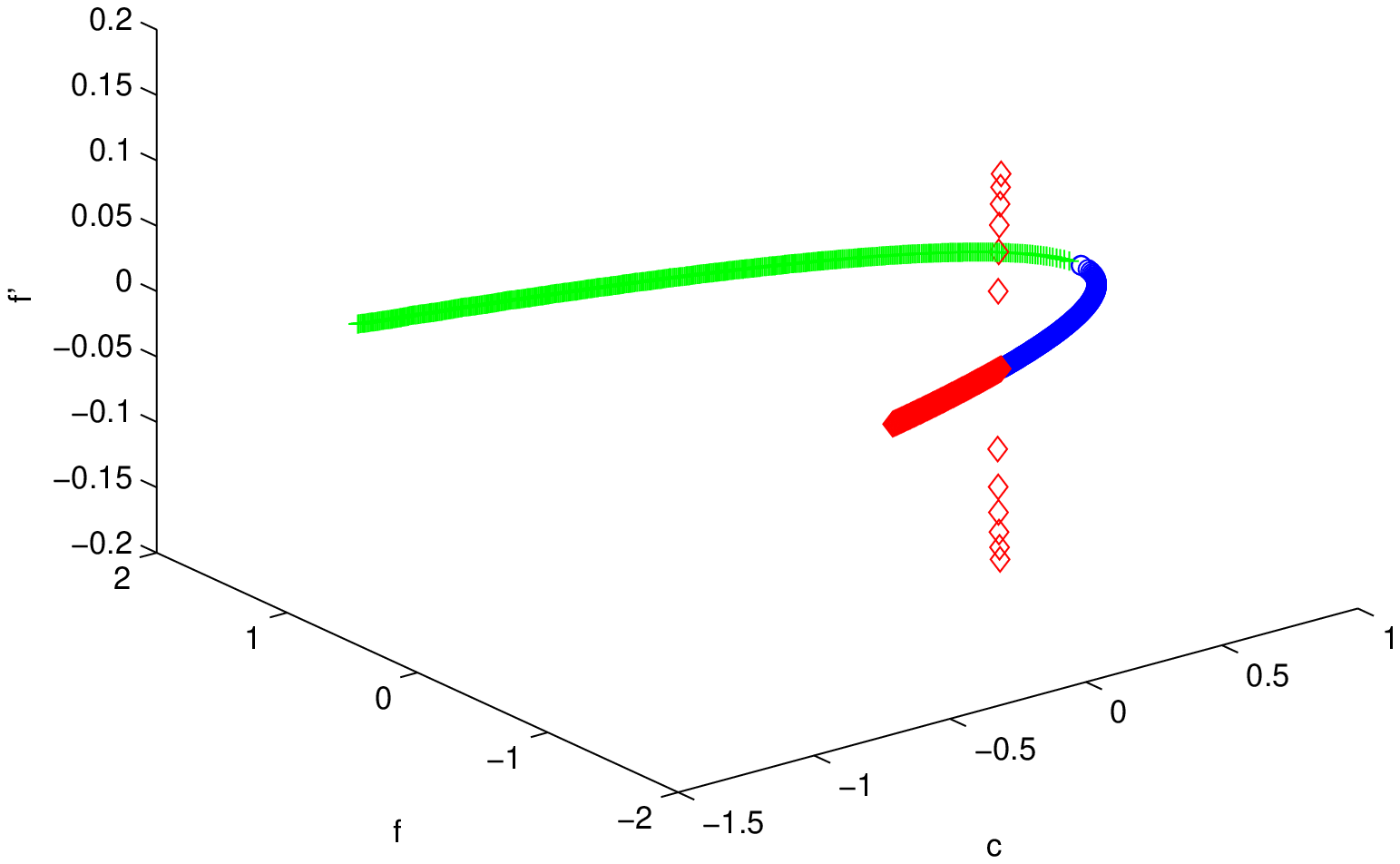}
\caption{Bifurcation diagram, coded by spectrum of
$\frac{d^2}{dx^2}-2f$: green = nonpositive spectrum, blue = one positive
eigenvalue, red = two positive eigenvalues}
\label{bif_diag}
\end{figure}

Once computed, the numerical solutions can then be tabulated
conveniently in a bifurcation diagram.  That is, consider the set in
$\mathbb{R}^3$ given by $(c,f(0),f'(0))$ for each solution $f$.
Evidently, by existence and uniqueness for ordinary differential
equations, each solution can be uniquely represented by such a point.
The results of such a computation are shown in Figure \ref{bif_diag}.
In this diagram, the solutions are color-coded by the number of
positive eigenvalues of $\frac{d^2}{dx^2}-2f$. \cite{Mazya_2005} (It
should be noted that the green curve continues for $c<-1.2$, but was
stopped for display reasons.)

\begin{figure}
\includegraphics[height=3in]{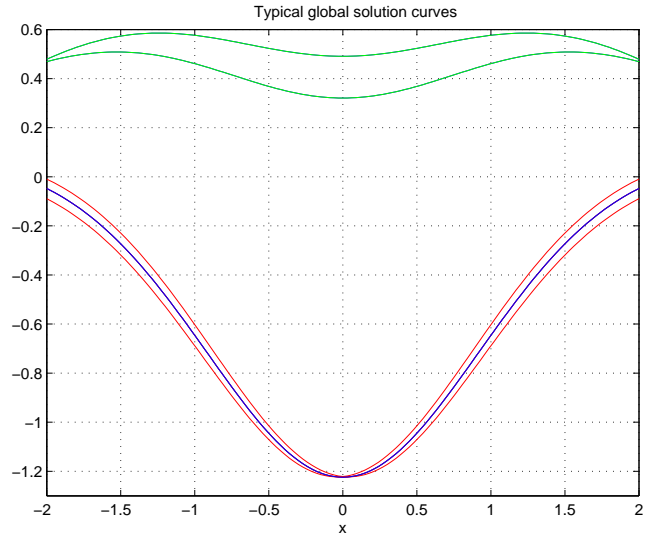}
\caption{Typical global solutions: green are from the positive branch,
  the blue one is taken from the lower branch with $f'(0)=0$, and the red
  ones are from the fork arms}
\end{figure}

Considering the bifurcation diagram, it appears to indicate that
\eqref{ode1} undergoes a saddle node bifurcation at approximately
$c=0.7706$, and a subcritical pitchfork bifurcation at $c=0.0501$.
The results agree with Theorem \ref{long_soln_exists}, in that
solutions do exist when $c<0$.  The saddle node bifurcation was
anticipated by the general shape of $\phi$.  For $c>0.7706$, global
solutions do not exist, which was qualitatively predicted
by Corollary \ref{phi_integral_bound}.

\begin{figure}
\includegraphics[height=3in]{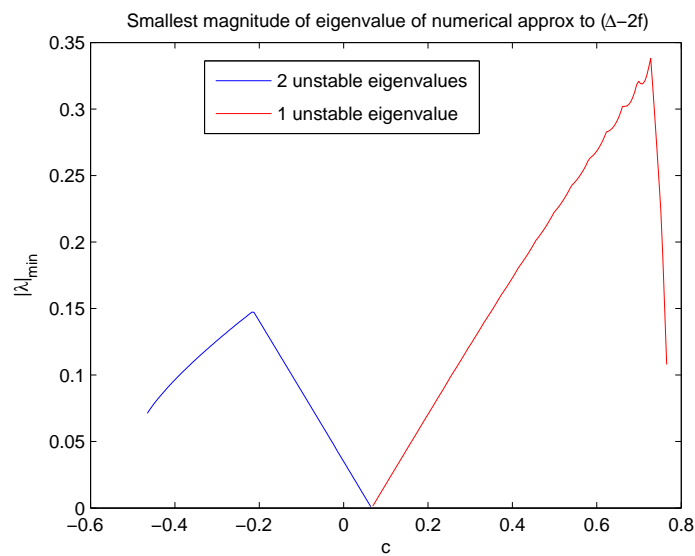}
\caption{Smallest-magnitude eigenvalue measured along the lower branch with $f'(0)=0$}
\label{small_eig}
\end{figure}

However, there are some stranger features of the bifurcation diagram.
Most prominently, the bifurcation diagram appears simply to {\it end}
near $c=-0.4652$, and at each branch of the pitchfork at $c=0.0740$.
It is important to verify that these are not numerical or
discretization errors.  If these {\it ends} are to be thought of as
valid bifurcations, very likely, $\frac{d^2}{dx^2} - 2f$ acquires a
zero eigenvalue there.  Plotting the smallest magnitude eigenvalue
gives some credence to this possibility.  (See Figure \ref{small_eig})

\begin{figure}
\includegraphics[height=3in]{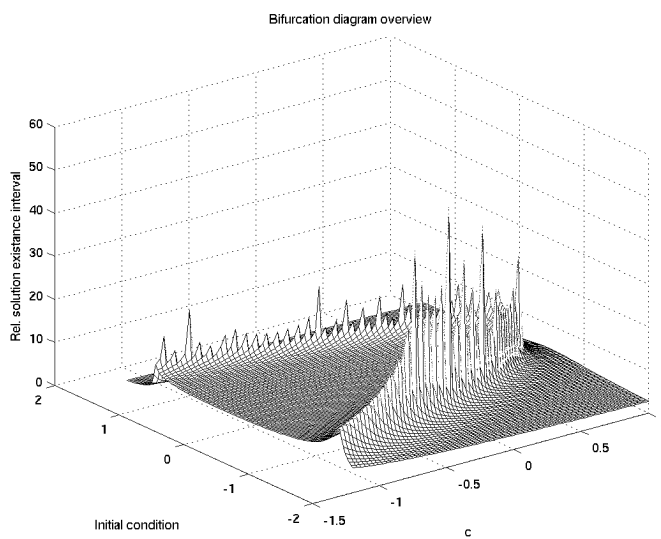}
\caption{Estimate of existence interval length}
\label{bif_exist}
\end{figure}

As another check, one can measure the size of the existence interval
for solutions to \eqref{long_ode1}, centered at $x=0$.  Looking in the
$(c,f(0))$-plane (taking $f'(0)=0$), one can find the first $x$ such
that the solution exceeds a particular value.  This is shown in Figure
\ref{bif_exist}, in which one sees the same general shape as in the
bifurcation diagram.  (The jagged nature of the graph along the actual
bifurcation diagram is due to aliasing.)  However, for $c<-0.4652$,
the lower branch clearly continues into solutions that exist for only
finite $x$.  So the end bifurcation indicates a failure of the
solutions to \eqref{long_ode1} to exist for all $x$.

\section{Conclusions}

In this article, an approach for counting and approximating global
solutions to a nonlinear, nonautonomous differential equation was
described that combines asymptotic and numerical information.  The
asymptotic information alone is enough to give necessary and
sufficient (but not sharp) conditions for solutions to exist, and
provides a fairly weak uniqueness condition.  More importantly, the
asymptotic approximation can be used to supply enough information to
pose a boundary value problem on a bounded interval containing a
smaller interval where asymptotic approximation is not valid.  This
boundary value problem is well-suited for numerical examination, and
the combined approach yields much more detailed results
than either method alone.

\section{Acknowledgements}
The author wishes to thank Dr. John H. Hubbard for suggesting this
kind of approach for approximating global solutions to ODE, and
Dr. Donna A. Dietz for implementing the numerical solver in MAPLE.

\bibliography{nonauto_bib}
\bibliographystyle{plain}

\end{document}